\documentclass[reqno,centertags, 12pt]{amsart}
\usepackage{amsmath,amsthm,amscd,amssymb}
\usepackage{latexsym}
\newcommand{\bbC}{{\mathbb{C}}}
\newcommand{\bbD}{{\mathbb{D}}}

\newcommand{\bbR}{{\mathbb{R}}}

\newcommand{\bbZ}{{\mathbb{Z}}}
\newcommand{\bdone}{{\boldsymbol{1}}}
\newcommand{\bddot}{{\boldsymbol{\cdot}}}

\newcommand{\calC}{{\mathcal{C}}}

\newcommand{\calG}{{\mathcal G}}
\newcommand{\calH}{{\mathcal H}}
\newcommand{\calL}{{\mathcal L}}
\newcommand{\calM}{{\mathcal M}}

\newcommand{\calR}{{\mathcal R}}

\newcommand{\lb}{\label}
\newcommand{\f}{\frac}

\newcommand{\ol}{\overline}
\newcommand{\ti}{\tilde  }
\newcommand{\wti}{\widetilde  }

\newcommand{\tr}{\text{\rm{Tr}}}
\newcommand{\dist}{\text{\rm{dist}}}
\newcommand{\loc}{\text{\rm{loc}}}

\newcommand{\spec}{\text{\rm{spec}}}

\newcommand{\ess}{\text{\rm{ess}}}

\newcommand{\bi}{\bibitem}

\newcommand{\beq}{\begin{equation}}
\newcommand{\eeq}{\end{equation}}
\newcommand{\ba}{\begin{align}}
\newcommand{\ea}{\end{align}}
\newcommand{\veps}{\varepsilon}

\newcounter{smalllist}
\newenvironment{SL}{\begin{list}{{\rm\roman{smalllist})}}{%
\setlength{\topsep}{0mm}\setlength{\parsep}{0mm}\setlength{\itemsep}{0mm}%
\setlength{\labelwidth}{2em}\setlength{\leftmargin}{2em}\usecounter{smalllist}%
}}{\end{list}}

\newcommand{\bigtimes}{\mathop{\mathchoice%
{\smash{\vcenter{\hbox{\LARGE$\times$}}}\vphantom{\prod}}%
{\smash{\vcenter{\hbox{\Large$\times$}}}\vphantom{\prod}}%
{\times}%
{\times}%
}\displaylimits}

\allowdisplaybreaks
\numberwithin{equation}{section}

\newtheorem{theorem}{Theorem}[section]

\newtheorem*{p2.1}{Proposition 2.1}
\newtheorem{proposition}[theorem]{Proposition}

\theoremstyle{definition}
\newtheorem{example}[theorem]{Example}

\theoremstyle{remark}
\newtheorem*{remark}{Remark} 
\newtheorem*{remarks}{Remarks}

\theoremstyle{definition}
\newtheorem*{definition}{Definition}

\newcommand{\abs}[1]{\lvert#1\rvert}

\begin{document}
\title[The Essential Spectrum]{The Essential Spectrum of Schr\"odinger, 
Jacobi, and CMV Operators}
\author[Y.~Last and B.~Simon]{Yoram Last$^{1,3}$ and Barry Simon$^{2,3}$}

\thanks{$^1$ Institute of Mathematics, The Hebrew University,
91904 Jerusalem, Israel. E-mail: ylast@math.huji.ac.il. Supported in part
by THE ISRAEL SCIENCE FOUNDATION (grant no. 188/02)}
\thanks{$^2$ Mathematics 253-37, California Institute of Technology, Pasadena, CA 91125.
E-mail: bsimon@caltech.edu. Supported in part by NSF grant DMS-0140592}  
\thanks{$^3$ Research supported in part 
by Grant No.\ 2002068 from the United States-Israel Binational Science Foundation 
(BSF), Jerusalem, Israel}

\date{April 12, 2005} 

\begin{abstract}  We provide a very general result that identifies the essential 
spectrum of broad classes of operators
as exactly equal to 
the closure of the union of the spectra of suitable limits at infinity. Included 
is a new result on the essential spectra when potentials are asymptotic to 
isospectral tori. We also recover with a unified framework the HVZ theorem 
and Krein's results on orthogonal polynomials with finite essential spectra. 
\end{abstract}

\maketitle

\section{Introduction} \lb{s1} 

One of the most simple but also most powerful ideas in spectral theory is Weyl's theorem, 
of which a typical application is (in this introduction, in order to avoid technicalities, 
we take potentials bounded):  

\begin{theorem}\lb{T1.1} If $V,W$\! are bounded functions on $\bbR^\nu$ and 
$\lim_{\abs{x}\to\infty}  [V(x)-W(x)] =0$, then 
\begin{equation} \lb{1.1} 
\sigma_\ess (-\Delta + V)=\sigma_\ess (-\Delta +W)  
\end{equation}  
\end{theorem} 

Our goal in this paper is to find a generalization of this result that allows ``slippage" near 
infinity. Typical of our results are the following: 

\begin{theorem}\lb{T1.2} Let $V$ be a bounded periodic function on $(-\infty,\infty)$ and 
$H_V$ the operator $-\f{d^2}{dx^2}+V(x)$ on $L^2 (\bbR)$. For $x>0$, define $W(x)=V(x + 
\sqrt{x})$ and let $H_W$ be $-\f{d^2}{dx^2}+W(x)$ on $L^2(0,\infty)$ with some selfadjoint 
boundary conditions at zero. Then 
\begin{equation} \lb{1.2} 
\sigma_\ess (H_W) = \sigma (H_V)  
\end{equation} 
\end{theorem} 

\begin{theorem} \lb{T1.3} Let $\alpha$ be irrational and let $H$ be the discrete Schr\"odinger 
operator on $\ell^2 (\bbZ)$ with potential $\lambda\cos(\alpha n)$. Let $\wti H$ be the 
discrete Schr\"odinger operator on $\ell^2 (\{0,1,2,\dots\})$ with potential $\lambda \cos 
(\alpha n + \sqrt{n})$. Then 
\begin{equation} \lb{1.3} 
\sigma_\ess (\wti H) = \sigma (H)  
\end{equation} 
\end{theorem} 

Our original motivation in this work was extending a theorem of Barrios-L\'opez \cite{BRLL}   
in the theory of orthogonal polynomials on the unit circle (OPUC); see \cite{OPUC1,OPUC2}.  

\begin{theorem}[see Example~4.3.10 of \cite{OPUC1}]\lb{T1.4A} Let $\{\alpha_n\}_{n=0}^\infty$ 
be a sequence of Verblunsky coefficients so that for some $a\in (0,1)$, one has
\begin{equation} \lb{1.4x} 
\lim_{n\to\infty}\, \abs{\alpha_n} =a \qquad 
\lim_{n\to\infty}  \f{\alpha_{n+1}}{\alpha_n}=1  
\end{equation}  
Then the CMV matrix for $\alpha_n$ has essential spectrum identical to the case  
$\alpha_n\equiv a$. 
\end{theorem} 

This goes beyond Weyl's theorem in that $\alpha_n$ may not approach $a$; rather $\abs{\alpha_n} 
\to a$ but the phase is slowly varying and may not have a limit. The way to understand this 
result is to realize that $\alpha_n\equiv a$ is a periodic set of Verblunsky coefficients. 
The set of periodic coefficients with the same essential spectrum is, for each $\lambda\in
\partial\bbD$ ($\bbD=\{z\mid\abs{z}<1\}$), the constant sequence $\alpha_n =\lambda a$. 
\eqref{1.4} says in a precise sense that the given $\alpha_n$ is approaching this  
isospectral torus. We wanted to prove, and have proven, the following: 

\begin{theorem}\lb{T1.4} If a set of Verblunsky coefficients or Jacobi 
parameters is asymptotic to an isospectral torus, then the essential spectrum 
of the corresponding CMV or Jacobi matrix is identical to the common essential 
spectrum of the isospectral torus. 
\end{theorem} 

In Section~\ref{s5}, we will be precise about what we mean by ``asymptotic to 
an isospectral torus." Theorem~\ref{T1.4} positively settles Conjecture~12.2.3 
of \cite{OPUC2}. 

In the end, we found an extremely general result. To describe it, we recall some 
ideas in our earlier paper \cite{LastS}. We will first consider Jacobi matrices 
($b_n\in\bbR$, $a_n>0$) 
\begin{equation} \lb{1.4} 
J=\begin{pmatrix} b_1 & a_1 & 0 & \cdots \\
a_1 & b_2 & a_2 & \cdots \\
0 & a_2 & b_3 & \ddots \\
\vdots & \vdots & \ddots & \ddots
\end{pmatrix}
\end{equation} 
where, in line with our convention to deal with the simplest cases in this introduction, we 
suppose there is a $K\in (0,\infty)$ so 
\begin{equation} \lb{1.5} 
\sup_n\, \abs{b_n} + \sup_n\, \abs{a_n} + \sup_n\, \abs{a_n}^{-1} \leq K  
\end{equation} 
A right limit point of $J$ is a double-sided Jacobi matrix, $J^{(r)}$, with parameters 
$\{a_n^{(r)},b_n^{(r)}\}_{n=-\infty}^\infty$ so that there is a subsequence $n_j$ with 
\begin{equation} \lb{1.6} 
a_{n_j +\ell} \to a_\ell^{(r)} \qquad 
b_{n_j+\ell} \to b_\ell^{(r)}  
\end{equation}  
as $j\to\infty$ for each fixed $\ell=0,\pm 1, \pm 2, \dots$. In \cite{LastS}, we noted that 

\begin{proposition}\lb{P1.5} For each right limit point, $\sigma (J^{(r)})\subset 
\sigma_\ess (J)$. 
\end{proposition} 

This is a basic result that many, including us, regard as immediate. For if 
$\lambda\in \sigma (J^{(r)})$ and $\varphi^{(m)}$ is a sequence of unit trial functions 
with $\|(J^{(r)}-\lambda)\varphi^{(m)}\|\to 0$, then for any $j(m)\to\infty$, 
$\|(J-\lambda) \varphi^{(m)} (\bddot + n_{j(m)})\|\to 0$, and if $j(m)$ is chosen 
going to infinity fast enough, then $\varphi^{(m)} (\bddot - n_{j(m)})\to 0$ weakly, 
so $\lambda\in\sigma_\ess (J)$. 

Let $\calR$ be the set of right limit points. Clearly, Proposition~\ref{P1.5} says that 
\begin{equation} \lb{1.7} 
\ol{\bigcup_{r\in\calR}\, \sigma (J^{(r)})} \subset \sigma_\ess (J)  
\end{equation} 
Our new realization here for this example is that 

\begin{theorem}\lb{T1.6} If \eqref{1.5} holds, then 
\begin{equation} \lb{1.8} 
\ol{\bigcup_{r\in\calR}\, \sigma (J^{(r)})} = \sigma_\ess (J)  
\end{equation} 
\end{theorem} 

\begin{remark}
It is an interesting question whether anything is gained in \eqref{1.8} by
taking the closure---that is, whether the union is already closed.  In every
example we can analyze the union is closed. V.~Georgescu has informed us that 
the methods of \cite{GI1} imply that the union is always closed and that the 
details of the proof of this fact are the object of a paper in preparation
\end{remark} 

Surprisingly, the proof will be a rather simple trial function argument. The difficulty 
with such an argument tried naively is the following: To say $J^{(r)}$ is a right limit 
point means that there are $L_m\to\infty$ so that $J\restriction [n_{j(m)} - L_m, 
n_{j(m)} +L_m]$ shifted to $[-L_m, L_m]$ converges uniformly to $J^{(r)}\restriction 
[-L_m,L_m]$. But $L_m$ might grow very slowly with $m$. Weyl's criterion says that if 
$\lambda\in\sigma_\ess (J)$, there are trial functions, $\varphi_k$, supported on 
$[n_k -\ti L_k, n_k + \ti L_k]$ so $\|(J-\lambda)\varphi_k\|\to 0$. By a compactness 
argument, one can suppose the $n_k$ are actually $n_{j(m)}$'s for some right limit. 
The difficulty is that $\ti L_m$ might grow much faster than $L_m$, so translated 
$\varphi_k$'s are not good trial functions for $J^{(r)}$. 

The key to overcoming this difficulty is to prove that one can localize trial functions 
in some interval of fixed size $L$, making a localization error of $O(L^{-1})$. 
This is what we will do in Section~\ref{s2}. In this idea, we were motivated by 
arguments in Avron et al.\ \cite{AvMS}, although to handle the continuum case, we 
will need to work harder. 

The use of localization ideas to understand essential spectrum, an implementation 
using double commutators, is not new---it goes back to Enss \cite{Enss} and 
was raised to high art by Sigal \cite{Sig}. Enss and Sigal, and also Agmon 
\cite{Agm} and Garding \cite{Gard}, later used these ideas and positivity inequalities 
to locate  $\inf\sigma_\ess(H)$, which suffices for the HVZ theorem but not for some 
of our applications. 

What distinguishes our approach and allows stronger results is that, first, we 
use trial functions exclusively and, second, as noted above, we study all of 
$\sigma_\ess$ rather than only its infimum. Third, and most significantly, we 
do not limit ourselves to sets that are cones near infinity and instead take balls. 
This gives us small operator errors rather than compact operator errors
(although one can modify 
our arguments and take ball sizes that go to infinity slowly, and so get a compact 
localization error). It makes the method much more flexible. 

While this paper is lengthy because of many different applications, the underlying 
idea is captured by the mantra ``localization plus compactness."  Here compactness 
means that resolvents restricted to balls of fixed size translated to zero lie in 
compact sets. We have in mind the topology of norm convergence once resolvents
are multiplied by the characteristic functions of arbitrary fixed balls.

Because we need to control $\|(A-\lambda)\varphi\|^2$ and not just $\langle 
\varphi, (A-\lambda)\varphi\rangle$, if we used double commutators, we would need 
to control $[j,[j,(A-\lambda)^2]]$, so in the continuum case we get unbounded operators and 
the double commutator is complicated. For this reason, following \cite{AvMS} and \cite{HK}, 
we use single commutators and settle for an inequality rather than the equality 
one gets from double commutators.  

After we completed this paper and released a preprint, we learned of some related 
work using $C^*$-algebra techniques to compute $\sigma_\ess(H)$ as the closure of 
a union of spectra of asymptotic Hamiltonians at infinity; see  
Georgescu-Iftimovici \cite{GI1} and Mantoiu \cite{M1}. Further work is in 
\cite{AMP,GG2,GI3,GI2,MPR,Rod}. 

We also learned of very recent work of Rabinovich \cite{Rab05}, based in part on 
\cite{Rab9,Rab6,Rab10,Rab5,Rab7,Rab8}, using the theory of Fredholm operators 
to obtain results on essential spectrum as a union of spectra of suitable limits 
at infinity. 

Thus, the notion that in great generality the essential spectra is a union of 
spectra of limits at infinity is not new. Our contributions are twofold: First, 
some may find our direct proof via trial functions more palatable than arguments 
relying on considerable machinery. Second, our examples of Section~\ref{snew4}, 
Section~\ref{s5}, and Section~\ref{s7}(b),(c) are, so far as we know, new, although 
it is certainly true that the methods of \cite{GI1,M1,Rab05}, can analyze some or 
even all these examples. In particular, we settle Conjecture~12.2.3 of \cite{OPUC2}. 

There is obviously considerable overlap in philosophies (which, after all, both extend 
the ideas of the HVZ theorem) and results. The techniques seem to be rather different, 
although we suspect a translation of the $C^*$-algebra machinery to more prosaic 
terms will show similarities that are, for now, not clear to us. 

The paper \cite{GI3} has results stated without reference to $C^*$-algebras (although 
the proofs use them) and, in particular, our Theorems~\ref{T3.7} and \ref{T3.11} are 
special cases of Theorem~1.1 of \cite{GI3}. 

We present the localization lemmas in Section~\ref{s2} and prove our main results 
in Section~\ref{s3}. Section~\ref{snew4} discusses an interesting phenomena involving 
Schr\"odinger operators with severe oscillations at infinity. Section~\ref{s5} has 
the applications to potentials asymptotic to isospectral tori and includes results 
stronger than Theorems~\ref{T1.2}, \ref{T1.3}, and \ref{T1.4}. In particular, we 
settle positively Conjecture~12.2.3 of \cite{OPUC2}. Section~\ref{s6} discusses the HVZ 
theorem, and Section~\ref{s7} other applications. Section~\ref{s8} discusses magnetic 
fields. 

We can handle the common Schr\"odinger operators associated to quantum theory with or
without magnetic fields as well as orthogonal polynomials on the real line (OPRL)
and unit circle (OPUC). 

\smallskip 
It is a pleasure to thank D.~Damanik and R.~Killip for useful discussions, and 
V.~Georgescu, M.~Mantoiu, V.~Rabinovich, A.~Sobolev and B.~Thaller for useful 
correspondence. This research was completed during B.~S.'s stay as a Lady Davis 
Visiting Professor at The Hebrew University of Jerusalem. He would like to thank 
H.~Farkas for the hospitality of the Einstein Institute of Mathematics at 
The Hebrew University. 

\section{Localization Estimates} \lb{s2}

Here we will use localization formulae but with partitions of unity that are concentrated 
on balls of fixed size in place of the previous applications that typically take $j$'s 
that are homogeneous of degree zero near infinity. Also, we use single commutators. 

Let $\calH$ be a separable Hilbert space and $A$ a selfadjoint operator on $\calH$.
Let $\{j_\alpha\}$ be a set of bounded selfadjoint operators indexed by either a
discrete set, $S$, like $\bbZ^\nu$ or by 
$\alpha\in\bbR^\nu$. In the latter case, we suppose $j_\alpha$ is measurable and 
uniformly bounded in $\alpha$.
We assume that $\{j_\alpha\}$ is a partition of unity, namely, 
\begin{equation} \lb{2.1} 
\sum_{\alpha\in S} j_\alpha^2 = \bdone  \qquad\text{or}\qquad 
\int_{\alpha\in\bbR^\nu} j_\alpha^2\, d^\nu\alpha =\bdone
\end{equation}
where the convergence of the sum or the meaning of the integral is in the
weak operator topology sense.  
Two examples that will often arise are where $\calH=\ell^2 (\bbZ^\nu)$, $\psi\in\ell^2 
(\bbZ^\nu)$ is real-valued with $\sum_n \psi(n)^2 =1$, and $\{j_m\}_{m\in\bbZ^\nu}$ is 
multiplication by $\psi (\bddot -m)$, or where $\calH=L^2 (\bbR^\nu,d^\nu x)$, 
$\psi\in L^2(\bbR^\nu, d^\nu x)\cap L^\infty (\bbR^\nu, d^\nu x)$ is real-valued 
with $\int \psi(x)^2\, d^\nu x = 1$, and $\{j_y\}_{y\in\bbR^\nu}$ is multiplication by 
$\psi(\bddot -y)$. 

Assume that for each $\alpha$, $j_\alpha$ maps the domain of $A$ to itself
and let $\varphi$ be a vector in the domain of $A$.
Notice that 
\begin{align} 
\|Aj_\alpha \varphi\|^2 &= \|(j_\alpha A+[A,j_\alpha])\varphi\|^2 \notag \\ 
&\leq 2 \|j_\alpha A\varphi\|^2 + 2\|[A,j_\alpha]\varphi\|^2 \lb{2.2} 
\end{align} 
Thus 

\begin{proposition} \lb{P2.1} 
\begin{equation} \lb{2.3} 
\sum_\alpha \, \|Aj_\alpha\varphi\|^2 \leq 2 \|A\varphi\|^2 + \langle \varphi, 
C\varphi  \rangle
\end{equation} 
where 
\begin{equation} \lb{2.4} 
C=2\sum_\alpha - [A,j_\alpha]^2  
\end{equation}  
\end{proposition} 

\begin{remark} Since $[A,j_\alpha]$ is skew-adjoint,
$-[A,j_\alpha]^2 = [j_\alpha,A]^* [j_\alpha,A] \geq 0$. 
\end{remark}

\begin{proof} \eqref{2.3} is immediate from \eqref{2.2} since 
\begin{equation} \lb{2.5} 
\sum_\alpha \, \|j_\alpha A\varphi\|^2 = \sum_\alpha \langle A\varphi, j_\alpha^2 
A\varphi\rangle = \|A\varphi\|^2  
\end{equation} 
and  
\begin{equation} \lb{2.6} 
\|[A,j_\alpha]\varphi\|^2 = -\langle \varphi, [A,j_\alpha]^2 \varphi\rangle   
\end{equation} 
\end{proof} 

\begin{theorem}\lb{T2.2} There exists an $\alpha$ so that $j_\alpha\varphi\neq 0$ and 
\begin{equation} \lb{2.7} 
\|Aj_\alpha\varphi\|^2 \leq \biggl\{ 2\biggl( \f{\|A\varphi\|}{\|\varphi\|}\biggr)^2 + 
\|C\|\biggr\} \|j_\alpha \varphi\|^2  
\end{equation} 
\end{theorem} 

\begin{proof} Call the quantity in $\{\,\,\}$ in \eqref{2.7} $d$. Then, since $\|\varphi\|^2 
=\sum_\alpha \|j_\alpha\varphi\|^2$, \eqref{2.3} implies 
\[
\sum_\alpha\, [\|Aj_\alpha\varphi\|^2 - d \|j_\alpha\varphi\|^2 ]\leq 0
\]
so at least one term with $\|j_\alpha\varphi\|\neq 0$ is nonpositive. 
\end{proof} 

To deal with unbounded $A$'s, we will want to suppose that $\sqrt{C}$ is $A$-bounded: 

\begin{theorem} Suppose $A$ is unbounded and 
\begin{equation} \lb{2.8} 
\langle \varphi,C\varphi\rangle \leq \delta (\|A\varphi\|^2 + \|\varphi\|^2)  
\end{equation} 
Then there is an $\alpha$ with $j_\alpha \varphi\neq 0$ so that 
\begin{equation} \lb{2.9} 
\|Aj_\alpha \varphi\|^2 \leq \biggl\{ (2+\delta) \,\f{\|A\varphi\|^2}{\|\varphi\|^2} 
+ \delta\biggr\} \|j_\alpha \varphi\|^2  
\end{equation}  
\end{theorem} 

\begin{proof} By \eqref{2.3} and \eqref{2.8}, we have 
\[
\sum_\alpha \, \|Aj_\alpha \varphi\|^2 \leq (2+\delta) \|A\varphi\|^2 + 
\delta \|\varphi\|^2 
\]
so, as before, \eqref{2.9} follows. 
\end{proof}

\section{The Essential Spectrum} \lb{s3}

This is the central part of this paper. We begin with Theorem~\ref{T1.6}, the simplest 
of the results: 

\begin{proof}[Proof of Theorem~\ref{T1.6}] We already proved \eqref{1.7} in the 
remarks after Proposition~\ref{P1.5}, so suppose $\lambda\in\sigma_\ess (J)$. Recall 
Weyl's criterion, $\lambda\in\sigma_\ess (J)\Leftrightarrow$ there exist unit vectors  
$\varphi_m\overset{w}{\longrightarrow} 0$ with $\|(J-\lambda)\varphi_m\|\to 0$. 

Given $\veps$, pick a trial sequence $\{\varphi_m\}$, such that each $\varphi_m$ is
supported in $\{n\mid n>m\}$, so that 
\begin{equation} \lb{3.1} 
\|(J-\lambda)\varphi_m\|^2 \leq \tfrac13\, \veps^2 \|\varphi_m\|^2  
\end{equation} 
which we can do, by Weyl's criterion, since $f_j \overset{w}{\longrightarrow} 0$ implies 
$\sum_{n<m} \abs{f_j(n)}^2\to 0$ for each $m$. 

For $L=1,2,3,\dots$, let 
\begin{equation} \lb{3.2} 
\psi_L(n) = \begin{cases} 
\f{n-1}{L} & n=1,2,\dots, L \\ 
\f{2L-1-n}{L} & n=L, L+1, \dots, 2L-1 \\
0 & n\geq 2L-1 
\end{cases}  
\end{equation} 
and let 
\begin{equation} \lb{3.3} 
c_L^2 =\sum_n \abs{\psi_L(n)}^2 
\end{equation} 
so that $c_L\sim L^{1/2}$ in the sense that for some $0 < a\leq b < \infty$,    
\begin{equation} \lb{3.3a} 
aL^{1/2} \leq c_L \leq b L^{1/2}  
\end{equation}
For $\alpha =1,2,\dots$, let 
\begin{equation} \lb{3.4} 
j_{\alpha,L}(n) = c_L^{-1} \psi_L (n+\alpha)  
\end{equation} 
so, by \eqref{3.3}, 
\begin{equation} \lb{3.5} 
\sum_\alpha j_{\alpha,L}^2 \equiv 1  
\end{equation}  

Since $\abs{\psi_L (n+1) -\psi_L(n)}\leq L^{-1}$, we see that 
\begin{equation} \lb{3.6} 
\abs{\langle\delta_n, [j_{\alpha,L},J]\delta_m\rangle} =  
\begin{cases} 
\sup_n \abs{a_n}\,\, c_L^{-1} L^{-1} &\text{if } \abs{n-m} =1\text{ and } 
\abs{n-\alpha-L} \leq L \\
0 &\text{otherwise} 
\end{cases}  
\end{equation} 
Therefore, $C\equiv \sum_\alpha 2[j_{\alpha,L},J]^2$ is a 5-diagonal matrix with 
matrix elements bounded by 
\begin{equation} \lb{3.7} 
2\cdot 2 (2L) c_L^{-2} L^{-2} \bigl( \, \sup_n\, \abs{a_n}\bigr)^2  
\end{equation} 
where the second two comes from the number of $k$'s that make a nonzero contribution to 
$\langle \delta_n, [j_{\alpha,L},J]\delta_k\rangle \langle\delta_k, 
[j_{\alpha,L},J]\delta_m\rangle$. By \eqref{3.3a}, there is a constant $K$ depending 
on $\sup_n\abs{a_n}$ so that 
\begin{equation} \lb{3.8x} 
\|C\|\leq KL^{-2}  
\end{equation}  
Picking $L$ so $KL^{-2} < \veps^2/3$, we see, by Theorem~\ref{T2.2}, there is a
$j_{\alpha_m}$ so $\|j_{\alpha_m}\varphi_m\|\neq 0$ and 
\begin{equation} \lb{3.8} 
\|(J-\lambda) j_{\alpha_m}\varphi_m\| \leq \veps \|j_{\alpha_m}\varphi_m\|  
\end{equation} 

The intervals 
\[
I_m = [\alpha_m +1, \alpha_m +2L-1] 
\]
which support $j_{\alpha_m}\varphi_m$, have fixed size, and move out to infinity since 
$I_m\subset \{n\mid n\geq m-L\}$. Since the set of real numbers with $\abs{b} + 
\abs{a} + \abs{a}^{-1} \leq K$ is compact and $L$ is finite, we can find a right 
limit point $J^{(r)}$ so that a subsequence of $J\restriction I_m$ translated by 
$\alpha_m +L$ converges to $J^{(r)}\restriction [1-L,L-1]$.
Using translations of the trial functions $j_{\alpha_m}\varphi_m$,
we find $\psi_m$ so 
\begin{equation} \lb{3.9} 
\lim_{m\to\infty}\, \f{\|(J^{(r)}-\lambda) \psi_m\|} 
{\|\psi_m\|} \leq \veps  
\end{equation} 
which means 
\begin{equation} \lb{3.10x} 
\dist (\lambda,\sigma (J^{(r)})) \leq \veps  
\end{equation}  
Since  $\veps$ is arbitrary, we have $\lambda\in \ol{\cup\sigma (J^{(r)})}$. 
\end{proof} 

We have been pedantically careful about the above proof so that below we can be 
much briefer and just relate to this idea as ``localization plus compactness" 
and not provide details. 

We turn next to the CMV matrices defined by a sequence of Verblunsky coefficients 
$\{\alpha_j\}_{j=0}^\infty$ with $\alpha_j\in\bbD$. We define the unitary $2\times 2$ 
matrix $\Theta(\alpha) = \left(\begin{smallmatrix} \ol{\alpha} & \rho \\ \rho & -\alpha 
\end{smallmatrix}\right)$ where $\rho = (1-\abs{\alpha}^2)^{1/2}$ and $\calL = 
\Theta_0\oplus \Theta_2\oplus \Theta_4 \oplus \cdots$, $\calM =\bdone\oplus\Theta_1 
\oplus\Theta_3+\cdots$, where $\bdone$ is a $1\times 1$ matrix and $\Theta_j=\Theta 
(\alpha_j)$. Then the CMV matrix is the unitary matrix $\calC=\calL\calM$. Given a 
two-sided sequence $\{\alpha_j\}_{j=-\infty}^\infty$, we define $\ti\calL = 
\cdots\Theta_{-2}\oplus\Theta_0\oplus\Theta_2$ and $\wti\calM = \Theta_{-1}\oplus 
\Theta_1\oplus\Theta_3\oplus\cdots$ on $\ell^2(\bbZ)$ where $\Theta_j$ acts on the 
span of $\delta_j$ and $\delta_{j+1}$. We set $\ti\calC=\ti\calL\wti\calM$. See 
\cite{OPUC1,OPUC2} for a discussion of the connection of CMV and extended CMV 
matrices to OPUC. 

In \cite{OPUC1,OPUC2}, $\ti\calC$ is used for the transpose of $\calC$ (alternate 
CMV matrix). Its use here is very different! 

If $\{\alpha_j\}_{j=0}^\infty$ is a set of Verblunsky coefficients with 
\begin{equation} \lb{3.10} 
\sup_j\, \abs{\alpha_j} <1  
\end{equation} 
we call $\{\beta_j\}_{j=-\infty}^\infty$ a right limit point if there is a sequence 
$m_j$ so that for $\ell=0,\pm 1, \dots$, 
\begin{equation} \lb{3.11} 
\lim_{j\to\infty}\, \alpha_{m_j+\ell} = \beta_\ell  
\end{equation} 
and we call $\ti\calC(\beta)$ a right limit of $\calC(\alpha)$. We have 

\begin{theorem}\lb{T3.1} Let $\calC(\alpha)$ be the CMV matrix of a sequence obeying 
\eqref{3.10}. Let $\calR$ be the set of right limit extended CMV matrices. Then 
\begin{equation} \lb{3.12} 
\sigma_\ess (\calC(\alpha)) = \ol{\bigcup_\calR\, \ti\calC (\beta)}  
\end{equation}  
\end{theorem} 

\begin{proof} The arguments of Section~\ref{s2} extend to unitary $A$ if $-[j_\alpha, A]^2$ 
is replaced by $[j_\alpha,A]^* [j_\alpha,A]$. Matrix elements of $[j_\alpha,\calC]$ are 
bounded by $\sup_{n,\abs{k}\leq 2} \abs{j_\alpha (n+k) - j_\alpha (n)}$ since 
$\calC$ has matrix elements bounded by $1$ and is 5-diagonal. Thus, $C$ is 9-diagonal, but 
otherwise the argument extends with no change since $\{\alpha\mid\abs{\alpha}\leq \sup_j 
\abs{\alpha_j}\}$ is a compact subset of $\bbD$. 
\end{proof} 

Next, we want to remove the condition that $\sup \abs{\alpha_j} <1$ in the OPUC case 
and the conditions $\sup\abs{b_j} <\infty$ and $\inf \abs{a_j}>0$ in the OPRL case. 
The key, of course, is to preserve compactness, that is, existence of limit points,  
and to do that, we need only extend the notion of right limit. 

If $\{\alpha_j\}_{j=-\infty}^\infty$ is a two-sided sequence in $\ol{\bbD}$, one can still 
define $\ti\calC (\alpha_j)$ since $\Theta (\alpha_j)$ makes sense. If $\abs{\alpha_j}=1$, 
then $\rho_j =0$ and $\Theta(\alpha_j)=\left(\begin{smallmatrix} \ol{\alpha}_j & 0 \\ 
0 & -\alpha_j \end{smallmatrix}\right)$ is a direct sum in such a way $\calL$ and $\calM$ 
both decouple into direct sums on $\ell^2 (-\infty, j]\oplus \ell^2 [j+1,\infty)$ so 
$\calC$ decouples. If a single $\alpha_j$ has $\abs{\alpha_j}=1$, we decouple into two 
semi-infinite matrices (both related by unitary transforms to ordinary CMV matrices), 
but if more than one $\alpha_j$ has $\abs{\alpha_j}=1$, there are finite direct summands. 

In any event, we can define $\ti\calC (\alpha_j)$ for $\{\alpha_j\}\in\bigtimes_{j=-\infty}^\infty 
\ol{\bbD}$ and define right limit points of $\calC(\alpha_j)$ even if $\sup\abs{\alpha_j} 
=1$. Since matrix elements of $\calC$ are still bounded by $1$, $\calC$ is still 5-diagonal 
and $\bigtimes_{j=-\infty}^\infty \ol{\bbD}$ is compact, we immediately have 

\begin{theorem} \lb{T3.2} With the extended notion of $\ti\calC$, Theorem~\ref{T3.1} holds 
even if \eqref{3.10} fails. 
\end{theorem} 

For bounded Jacobi matrices, we still want $\sup (\abs{a_n} + \abs{b_n})<\infty$, but we 
do not need $\inf \abs{a_n}>0$. Again, the key is to allow two-sided Jacobi matrices, $J_r$, 
with some $a_n=0$, in which case $J_r$ decouples on $\ell^2 (-\infty,n]\oplus \ell^2[n+1,\infty)$. 
If a single $a_n=0$, there are two semi-infinite matrices.
If more than one $a_n=0$, there are finite 
Jacobi summands. Again, with no change in proof except for the change in the meaning 
of right limits to allow some $a_n^{(r)}=0$, we have 

\begin{theorem}\lb{T3.3} Theorem~\ref{T1.6} remains true if \eqref{1.5} is replaced by 
\begin{equation} \lb{3.13} 
\sup_n\, (\abs{a_n} + \abs{b_n}) <\infty
\end{equation} 
so long as $J^{(r)}$ are allowed with some $a_n^{(r)}=0$. 
\end{theorem} 

In Section~\ref{s7}, we will use Theorems~\ref{T3.2} and \ref{T3.3} to complement the 
analysis of Krein (which appeared in Akhiezer-Krein \cite{AK}) for bounded Jacobi matrices 
with finite essential spectrum, and of Golinskii \cite{Gol2000} for OPUC with finite 
derived sets. 

Our commutator argument requires that $\abs{a_n}$ is bounded, but one can also handle 
$\limsup \abs{b_n} =\infty$. It is useful to define: 

\begin{definition} Let $A$ be a possibly unbounded selfadjoint operator. We say that $+\infty$ 
lies in $\sigma_\ess (A)$ if $\sigma(A)$ is not bounded above, and $-\infty$ lies in 
$\sigma_\ess(A)$ if $\sigma(A)$ is not bounded below. 
\end{definition} 

We now allow two-sided Jacobi matrices, $\ti J$, with $b_n=+\infty$ and/or $b_n=-\infty$ 
(and also $a_n=0$). If $\abs{b_n}=\infty$, we decouple into $\ell^2 (-\infty, n-1]\oplus 
\ell^2 [n+1,\infty)$ and place $b_n$ in ``$\sigma_\ess (\ti J)$." With this extended definition, 
we still have compactness, that is, for any intervals in $\bbZ_+$, $I_1, I_2,\dots$ of fixed 
finite size, $\ell$, with $\ell^{-1} \sum_{j\in I_n} j\to\infty$, there is a subsequence 
converging to a set of Jacobi parameters with possibly $b_n =+\infty$ or $b_n=-\infty$. We 
therefore have 

\begin{theorem}\lb{T3.4} Theorem~\ref{T1.6} remains true if \eqref{1.5} is replaced by 
\begin{equation} \lb{3.14} 
\sup_n\, \abs{a_n} <\infty  
\end{equation} 
so long as $J^{(r)}$ are allowed to have some $a_n^{(r)}=0$ and/or some $b_n^{(r)}=\pm\infty$. 
\end{theorem} 

\begin{remarks} 1. This includes the conventions on when $\pm\infty$ lies in $\sigma_\ess(J)$. 
To prove this requires a simple separate argument. Namely, $\langle\delta_n, J\delta_n\rangle 
=b_n$, so $b_n\in$ numerical range of $J=$ convex hull of $\sigma(J)$. Thus, if $b_{n_j}
\to\pm\infty$, then $\pm\infty\in\sigma (J)$. 

2. If $\sup_n \abs{a_n}=\infty$, $\sigma_\ess$ can be very subtle; see \cite{JN02,JNS04}. 
\end{remarks} 

Next, we turn to Jacobi matrices on $\bbZ^\nu$ (including $\nu=1$), that is, $J$ acts 
on $\ell^2 (\bbZ^\nu)$ via 
\begin{equation} \lb{3.15} 
(Ju)(n) =\sum_{\abs{m-n}=1}\, a_{(n,m)} u(m) + \sum_n b_n u(n)  
\end{equation} 
where the $b_n$'s are indexed by $n\in\bbZ^\nu$ and the $a_{(n,m)}$'s by bonds
$\{m,n\}$ (unordered pairs) with $\abs{m-n}=1$. For simplicity of exposition, we suppose 
\begin{equation} \lb{3.16} 
\sup_{\abs{m-n}=1}\, (\abs{a_{(n,m)}}+\abs{a_{(n,m)}}^{-1}) + \sup_n\, \abs{b_n} <\infty   
\end{equation} 
although we can, as above, also handle some limits with $a_{(n,m)}=0$ or some $\abs{b_n}=\infty$. 
With no change, one can also control finite-range off-diagonal terms, and with some effort 
on controlling $[j_\alpha,J]$, it should be possible to control infinite-range off-diagonal 
terms with sufficiently rapid off-diagonal decay. 

Let us call $\ti J$ a limit point of $J$ at infinity if and only if there are points 
$n_j\in\bbZ^\nu$ with $n_j\to\infty$ so that for every finite $k,\ell$, 
\begin{equation} \lb{3.17} 
b_{n_j+\ell}\to \ti b_\ell \qquad 
a_{(n_j+\ell, n_j +k)} \to \ti a_{(k,\ell)}
\end{equation} 
Let $\calL$ denote the set of limits $\ti J$. Then 

\begin{theorem}\lb{T3.5}
Let $J$ be a Jacobi matrix of the form \eqref{3.15} on $\ell^2 (\bbZ^\nu)$. 
Suppose \eqref{3.16} holds. Then 
\begin{equation} \lb{3.18} 
\sigma_\ess (J) = \ol{\bigcup_{\ti J\in\calL}\, \sigma (\ti J)}
\end{equation}  
\end{theorem} 

\begin{proof} We can define partitions of unity $j_{\alpha,L}$ indexed by $\alpha\in\bbZ^\nu$ 
with $j_\alpha (n)\neq 0$ only if $\abs{n-\alpha}\leq L$ and with $-\sum_\alpha [j_\alpha,J]^2$ 
bounded by $O(L^{-2})$. With this, the proof is the same as in the one-dimensional case. 
\end{proof} 

It is often comforting to only consider limit points in a single direction. Because the 
sphere is compact, this is easy. 

\begin{definition} Let $e\in S^{\nu -1}$, the unit sphere in $\bbR^\nu$. We say $\ti J$ 
is a limit point in direction $e$ if the $n_j$ in \eqref{3.17} obey $n_j/\abs{n_j}\to e$. 
We let $\calL_e$ denote the limit points in direction $e$. 
\end{definition}

Suppose $\ti J$ is a limit point for $J$ with sequence $n_j$. Since $S^{\nu-1}$ 
is compact, we can find a subsequence $n_{j(k)}$ so $n_{j(k)}/\abs{n_{j(k)}}\to 
e_0$ for some $e_0$. The subsequence also converges to $\ti J$ so $\ti J$ is a limit 
point for direction $e_0$. Thus,  

\begin{theorem}\lb{T3.6}
Let $J$ be a Jacobi matrix of the form \eqref{3.15} on $\ell^2 (\bbR^\nu)$. 
Suppose \eqref{3.16} holds. Then 
\begin{equation} \lb{3.19} 
\sigma_\ess (J) = \ol{\bigcup_{e\in S^{\nu -1}}\, \bigcup_{\ti J\in\calL_e}\, \sigma(\ti J)}  
\end{equation} 
\end{theorem} 

For example, if $\nu=1$, we can consider left and right limit points. 

Finally, we turn to Schr\"odinger operators. Here we need some kind of compactness 
condition of the $-\Delta +V$ that prevents $V$ from oscillating wildly at infinity 
(but see the next section). We begin with a warmup case that will be the core of 
our general case: 

\begin{theorem}\lb{T3.7} Let $V$ be a uniformly continuous, bounded function on $\bbR^\nu$. 
For each $e\in S^\nu$, call $W$\! a limit of $V$ in direction $e$ if and only if there 
exists $x_j\in\bbR^\nu$ with $\abs{x_j}\to\infty$ and $x_j/\abs{x_j}\to e$ so that  
$V(x_j +y)\to W(y)$. Then, with $\calL_e$ the limits in direction $e$, 
\begin{equation} \lb{3.20} 
\sigma_\ess (-\Delta +V) = \ol{\bigcup_e\, \bigcup_{W\in\calL_e}\, \sigma(-\Delta +W)}
\end{equation} 
\end{theorem} 

\begin{remarks} 1. While we have not stated it explicitly, there is a result for 
half-line operators. 

2. Uniform continuity means $\forall \veps, \exists \delta$, so $\abs{x-y}< 
\delta \Rightarrow \abs{V(x)-V(y)} <\veps$. It is not hard to see this is 
equivalent to $\{V (\bddot +y)\}_{y\in\bbZ^\nu}$ being equicontinuous.  

3. This result appears in a more abstract formulation in Georgescu-Iftimovici 
\cite{GI3}. 
\end{remarks} 

\begin{proof} As noted, uniform continuity implies uniform equicontinuity so, 
by the Arzela-Ascoli theorem (see \cite{RS1}), given any sequence of balls $\{x\mid 
\abs{x-y_j} \leq L\}$, there is an $e$ and a $W$ in $\calL_e$ so $V(\bddot + y_j) 
\to W(\bddot)$ uniformly on $\{x\mid\abs{x}\leq L\}$. This is the compactness 
needed for our argument. 

To handle localization, pick any nonnegative rotation invariant
$C^\infty$ function $\psi$ supported on 
$\{x\mid\abs{x} \leq 1\}$ with $\int \psi(x)^2 \, d^\nu x=1$.
Define $j_{x,L}$ as the operator of multiplication by the function
\[ 
j_{x,L}(y) = L^{-\nu/2} \psi (L^{-1} (y-x)) 
\]
and note that 
\[ 
\int j_{x,L}^2\, d^\nu x =1 
\]
With $A=(-\Delta +V-\lambda)$ and $C=2\int -[A,j_{x,L}]^2\, d^\nu x$, we have 
\eqref{2.8} with $\delta =O(L^{-2})$, since
$C=L^{-2}(c_1\Delta+c_2)$ for constants $c_1$ and $c_2$
(for $C$ is translation and rotation invariant and scale
covariant).

\eqref{3.20} follows in the usual way. 
\end{proof}

Our final result in this section concerns Schr\"odinger operators with potentially 
singular $V$'s. As in the last case, we will suppose regularity at infinity. In the 
next section, we will show how to deal with irregular oscillations near infinity. 
Recall the Kato class and norm \cite{S135,CFKS} is defined by 

\begin{definition} $V\colon\bbR^\nu\to\bbR$ is said to live in the Kato class, 
$K_\nu$, if and only if 
\begin{equation} \lb{3.22} 
\lim_{\alpha\downarrow 0}\, \biggl[\, \sup_x \int_{\abs{x-y} \leq \alpha} 
\abs{x-y}^{2-\nu} \abs{V(y)}\, d^\nu y\biggr] =0  
\end{equation} 
(If $\nu=1,2$, the definition is different. If $\nu=2$, $\abs{x-y}^{2-\nu}$ 
is replaced by $\log [\abs{x-y}^{-1}]$, and if $\nu=1$, we require $\sup_x 
\int_{\abs{x-y}\leq 1} \abs{V(y)}\, dy<\infty$.) The $K_\nu$ norm is defined by 
\begin{equation} \lb{3.23} 
\|V\|_{K_\nu} = \sup_x \int_{\abs{x-y}\leq 1} \abs{x-y}^{2-\nu} \abs{V(y)}\, d^\nu y  
\end{equation}  
\end{definition} 

We introduce here 

\begin{definition} $V\colon\bbR^\nu\to\bbR$ is called uniformly Kato if and only if 
$V\in K_\nu$ and 
\begin{equation} \lb{3.24} 
\lim_{y\downarrow 0}\, \|V-V(\bddot -y)\|_{K_\nu}=0  
\end{equation} 
\end{definition} 

\begin{example}\lb{E3.8} Let 
\begin{equation} \lb{3.25} 
V(x) =\sin (x_1^2)  
\end{equation} 
Then $V\in K_\nu$, but for $(x_0)_1$ large and $y=(\pi/2(x_0)_1, y_2, \dots)$, 
$[(x_0 +y)_1]^2 = (x_0)_1^2 + \pi + O(1/(x_0)_1)$, so for $x$ near $x_0$, $V(x)-
V(x-y)\sim 2V(x)$, and because of the $\abs{V(\bddot)}$ in \eqref{3.23}, we do not 
have \eqref{3.24}. We discuss this further in the next section. 
\end{example} 

\begin{example} \lb{E3.9} We say $p$ is canonical for $\bbR^\mu$ if $p=\mu/2$ 
where $\mu \geq 3$, $p>2$ if $\mu =2$, and $p=1$ if $\mu =1$. If 
\begin{equation} \lb{3.26} 
\sup_x \int_{\abs{x-y} \leq 1} \abs{V(y)}^p \, d^\mu y <\infty
\end{equation}  
then $V\in K_\mu$ (see \cite{CFKS}). Moreover, if 
\begin{equation} \lb{3.27} 
\lim_{\abs{x}\to\infty} \int_{\abs{x-y}\leq 1} \abs{V(y)}^p \, d^\mu y =0
\end{equation} 
it is easy to see that \eqref{3.24} holds because $V$ is small at infinity, 
and \eqref{3.24} holds for $L^p$ norm if $V\in L^p$. 
\end{example} 

\begin{example} \lb{E3.10} If $\pi\colon\bbR^\nu\to\bbR^\mu$ is a linear map onto 
$\bbR^\mu$ and $W\in K_\mu$, then $V(x)=W(\pi x)$ is in $K_\nu$ and the $K_\nu$ norm 
of $V$ is bounded by a $\pi$-dependent constant times the $K_\mu$ norm of $W$\!. 
If $W$ obeys \eqref{3.24}, so does $V$\!. 
\end{example}

We will combine Examples~\ref{E3.9} and \ref{E3.10} in our study of the HVZ theorem. 

\begin{proposition} \lb{P3.10} Let $V$ be a uniformly Kato potential on $\bbR^\nu$ 
and let $H_x =-\Delta + (\bddot -x)$. Then for any sequence $x_k\to\infty$, 
there is a subsequence $x_{k(m)}$ and a selfadjoint operator $H_\infty$ so that for 
$z\in\bbC\backslash [a,\infty)$ for some $a\in\bbR$, we have 
\begin{equation} \lb{3.28} 
\bigl\| [(H_{x_{k(m)}} -z)^{-1} -(H_\infty -z)^{-1}] \chi_S\|\to 0  
\end{equation} 
for $\chi_S$, the characteristic function of an arbitrary bounded set. 
\end{proposition} 

\begin{remark} Formally, $H_\infty$ is a Schr\"odinger operator of the form 
$H_0 +V_\infty$, but $V_\infty$, as constructed, is only in the completion of 
$K_\nu$, and that is known to include some distributions (see \cite{Gul,MV}). 
\end{remark} 

\begin{proof} It is known that if $W\in K_\nu$, then $W$ is $-\Delta$ form bounded 
with relative bound zero with bounds depending only on $K_\nu$ norms (see \cite{CFKS}). 
Thus, since all $V_x$'s have the same $K_\nu$ norm, we can find $a$ so $H_x\geq a$ for 
all $x$. It also means that for each $z\in\bbC\backslash [a,\infty)$, we can bound 
$\|\abs{W}^{1/2} (H_x-z)^{-1}\Delta^{1/2}\|$ by  $c \|W\|_{K_\nu}$ with $c$ only 
$z$-dependent and $\|V\|_{K_\nu}$-dependent. 

Let $\varphi$ be a $C^\infty$ function of compact support and note (constants are 
$z$- or $\|V\|_{K_\nu}$-dependent)  
\begin{align*}
\|\abs{W}^{1/2} [(H-z)^{-1}, \varphi]\| &\leq \|W^{1/2} (H-z)^{-1} [\Delta,\varphi] 
(H-z)^{-1}\|  \\
&\leq c \|W^{1/2} (H-z)^{-1} \Delta^{1/2}\|\, \|\nabla\varphi\| 
\end{align*} 
This in turn implies that if $S_1$ is a ball of radius $r$ fixed about $x_0$ and 
$S_2$ a ball of radius $R>r$, then 
\[
\|W^{1/2} (1-\chi_{S_2})(H-z)^{-1} \chi_{S_1}\| \to 0 
\]
as $R\to\infty$. So if $\|(W_n -W)\chi_S\|_{K_\nu}\to 0$ for all balls, and $\sup_n 
\|W_n\|_{K_\nu}<\infty$, then 
\[
\|((-\Delta +W_n-z)^{-1} - (-\Delta +W-z)^{-1}) \chi_S\| \to 0 
\]
for all $S$. 

In this way, we see that if $V$ is uniformly Kato and $V_{x_n}\to V_\infty$ 
in $K_\nu$ uniformly on all balls, then 
\begin{equation} \lb{3.29} 
\| [(H_{x_n}-z)^{-1} - (H_\infty -z)]^{-1} \chi_S\|\to 0  
\end{equation}  

The condition of $V$ being uniformly Kato means convolutions of $V$ with a $C^\infty$ 
approximate identity converge to $V$ in $K_\nu$ norm. Call the approximations $V^{(m)}$. 
Each is $C^\infty$ with bounded derivatives and so, by the equicontinuity argument in 
Theorem~\ref{T3.7}, we can find $x_{j_m (n)}$ and $V_\infty^{(m)}$ so 
\[ 
\| [ (-\Delta + V_{x_{j_m(n)}}^{(m)}-z)^{-1} - (-\Delta + V_\infty^{(m)}-z)^{-1}] 
\chi_S\| \to 0
\]
Since $V_x^{(m)}\to V_x$ uniformly in $x$, a standard $\veps/3$ argument (see \cite{RS1}) 
shows that one can find $x_{j(m)}$ so $\|[(H_{x_m}-z)^{-1} - (H_{x_{m'}}-z)^{-1}] 
\chi_S\|$ is small for each $S$ as $m,m'\to \infty$. In this way, we obtain the necessary 
limit operator. 
\end{proof} 

Given $V$ uniformly Kato, the limits constructed by Proposition~\ref{P3.10} where 
$x_n/\abs{x_n}\to e$ are called limits of $H$ in direction $e$. Again, the next result 
appears in a more abstract setting in Georgescu-Iftimovici \cite{GI3}.

\begin{theorem}\lb{T3.11} Let $V$ be uniformly Kato. Let $\calL_e$ denote the limits 
of $H$ in direction $e$. Then, 
\begin{equation} \lb{3.30} 
\sigma_\ess (H) = \ol{\bigcup_e\, \bigcup_{H_\infty\in\calL_e} \sigma(H_\infty)} 
\end{equation} 
\end{theorem} 

The papers that use $C^*$-algebras \cite{GI1,M1} study $h(p)+V$ in place of 
$-\Delta +V$. These papers only required that $h(p)\to\infty$ as $p\to\infty$. It seems  
likely that for many such $h$'s, our methods will work. $h(p)\to\infty$ as $p\to\infty$ 
is critical in our approach to assure that if $\varphi_n\to 0$ weakly and $\|(H-E) 
\varphi_n\|\to 0$ then $\chi_{\{x\mid\abs{x}<R\}}\varphi_n\to 0$. 

It is likely that one can develop a theory for $h(p) +V(x)$ without supposing $h(p) 
\to\infty$ or even $f(p,x)$, but one would need to consider limits at infinity in 
phase space, not just on configuration space. 

\begin{proof} We pick $a$ so $H_x\geq a$ for all $x$. Pick $z\in (-\infty,a)$ and 
let $\ti A_x =(H_x -z)^{-1}$. As above, $\|[\ti A_x,j_\alpha]\|  \leq c\|\nabla j_\alpha\|$ 
for any $j_\alpha$ in $C_0^\infty$. For $\lambda\in\sigma_\ess (H)$, let $A=(H_x -z)^{-1} 
-(\lambda-z)^{-1}$. Theorem~\ref{T2.2} provides the necessary localization estimate. 
Proposition~\ref{P3.10} provides the necessary compactness. \eqref{3.30} is then proven 
in the same way as earlier theorems. 
\end{proof}

\section{Schr\"odinger Operators With Severe Oscillations at Infinity} \lb{snew4}

This section is an aside to note that the lack of uniformity at infinity that can occur if 
$V$ is merely $K_\nu$ is irrelevant to essential spectrum. We begin with Example~\ref{E3.8}, 
the canonical example of severe oscillations at infinity: 

\begin{proposition} \lb{P4.1} Let 
\begin{equation} \lb{4.1} 
W(x)=\sin(x^2)  
\end{equation} 
on $(0,\infty)$ and let $H_0=-\f{d^2}{dx^2}$ with $u(0)=0$ boundary conditions. Then 
\begin{SL} 
\item[{\rm{(1)}}] $W(H_0+1)^{-1}$ is not compact. 
\item[{\rm{(2)}}] $(H_0 +1)^{-1/2} W(H_0-1)^{-1/2}$ is compact. 
\end{SL} 
\end{proposition} 

\begin{remarks} 1. Our proof of (1) shows that $Wf(H_0)$ is noncompact for any continuous 
$f\not\equiv 0$ on $(0,\infty)$. 

2. Consideration of $W=\vec{\nabla}\,\bddot\,\vec Q$ potentials goes back to the 1970's; 
(see \cite{BC,Cha,CM,Com,CG,DHS,Ism,IM,MS,Sar,Sch,Skr}). 
\end{remarks} 

\begin{proof} (1) Let $\varphi$ be a nonzero $C_0^\infty (0,\infty)$ function in $L^2$ and let 
\begin{equation} \lb{4.2} 
\psi_n(x) = [(H_0 +1) \varphi] (x-n)  
\end{equation} 
Then 
\begin{align} 
\|W(H_0 +1)^{-1}\psi_n \|^2 &= \int W(x)^2 \varphi(x-n)^2\, dx \notag \\
&= \tfrac12 \int \varphi(x)^2\, dx - \tfrac12 \int \cos(2x^2) \varphi(x-n)^2\, dx \notag \\ 
&\to \tfrac12 \int \varphi(x)^2\, dx \neq 0 \lb{4.3} 
\end{align} 
by an integration by parts. Since $\psi_n \overset{w}{\longrightarrow} 0$, this shows 
$W(H_0-1)^{-1}$ is not compact. 

\smallskip 
(2) Since $\f{d}{dx} [-\f{1}{2x} \cos (x^2)]=\sin(x^2) + O(x^{-2})$, we see $Q(x) = 
\lim_{y\to\infty} -\int_x^y W(z)\, dz$ exists and obeys 
\begin{equation} \lb{4.4} 
\abs{Q(x)}\leq c(x+1)^{-1}  
\end{equation} 
Thus $W=[\f{d}{dx},Q]$, so 
\[
(H_0+1)^{-1/2} W(H_0+1)^{-1/2} = \biggl((H_0+1)^{-1/2} \, \f{d}{dx}\biggr) (Q(H_0+1)^{-1/2}) + cc 
\]
Since $(H_0+1)^{-1/2} \f{d}{dx}$ is bounded and $Q(H_0-1)^{-1/2}$ is compact (by \eqref{4.4}), 
$(H_0+1)^{-1/2} W (H_0+1)^{-1/2}$ is compact. 
\end{proof} 

Thus, oscillations at infinity are irrelevant for essential spectrum! While the slick argument 
above somewhat obscures the underlying physics, the reason such oscillations do not matter 
has to do with the fact that $\sigma_\ess (H)$ involves fixed energy, and oscillations only 
matter at high energy. Our proof below will implement this strategy more directly. 

We begin by noting that the proof of Proposition~\ref{P3.10} implies the following: 

\begin{theorem}\lb{T4.2} Suppose $V_n$ is a sequence of multiplicative operators so that 
\begin{SL} 
\item[{\rm{(i)}}] For any $\veps >0$, there is $C_\veps$ so that 
\begin{equation} \lb{4.5} 
\langle\varphi, \abs{V_n}\varphi\rangle \leq \veps \|\nabla\varphi\|^2 + C_\veps \|\varphi\|^2  
\end{equation}  
for any $n$ and all $\varphi\in Q(-\Delta)$. 

\item[{\rm{(ii)}}] For any ball $S$ about zero, 
\begin{equation} \lb{4.6} 
\|(-\Delta +1)^{-1/2} \chi_S (V_n-V_m) (-\Delta +1)^{-1/2}\| \to 0  
\end{equation} 
as $n,m\to\infty$. 
\end{SL}
Then for any ball and $z\in\bbC\backslash [a,\infty)$, 
\begin{equation} \lb{4.7} 
\|[(-\Delta + V_n -z)^{-1} -(-\Delta +V_m-z)^{-1}]\chi_S\|\to 0
\end{equation} 
Moreover, if \eqref{4.6} holds as $n\to\infty$ with $V_m$ replaced by some $V_\infty$, then 
\begin{equation} \lb{4.8} 
\lim_{n\to\infty} \, \|[(-\Delta +V_n-z)^{-1} -(-\Delta +V_\infty -z)^{-1}] \chi_S\|=0  
\end{equation} 
\end{theorem} 

As an immediate corollary, we obtain 

\begin{theorem} \lb{T4.3} Let $V\in K_\nu$ obey 
\begin{equation} \lb{4.9} 
\lim_{R\to\infty} \, \sup_{\abs{x}\geq R}\, \int_{\abs{x-y}\leq 1} \abs{x-y}^{-(\nu-2)} 
\abs{V(y)}\, d^\nu y =0
\end{equation}  
Then 
\begin{equation} \lb{4.10} 
\sigma_\ess (-\Delta+V) = [0,\infty) 
\end{equation} 
\end{theorem} 

\begin{remark} If \eqref{4.9} holds, we say that $V$ is $K_\nu$ small at infinity. 
\end{remark} 

\begin{proof} By Theorem~\ref{T4.2}, if $x_n\to\infty$, 
\begin{equation} \lb{4.11} 
\|[(-\Delta +V (\bddot -x_n)-z)^{-1} -(-\Delta -z)^{-1}]\chi_S\|\to 0
\end{equation} 
so, in a sense, $-\Delta$ is the unique limit point at infinity. The standard localization 
argument proves \eqref{4.10}. 
\end{proof} 

Here is the key to studying general $V\in K_\nu$ with no uniformity at infinity: 

\begin{proposition} \lb{P4.4} Let $V_n$ be a sequence of functions supported in a fixed 
ball $\{x\mid\abs{x}\leq R\}$. Suppose 
\begin{equation} \lb{4.12} 
\lim_{\alpha\downarrow 0}\, \sup_{n,x}\, \int_{\abs{x-y}\leq\alpha} \abs{x-y}^{-(\nu -2)} 
\abs{V_n(y)}\, d^\nu y =0  
\end{equation} 
Then there is a subsequence $V_{n(j)}$ so 
\begin{equation} \lb{4.13} 
\lim_{j,k\to\infty} \, \|(-\Delta +1)^{-1/2} (V_{n(j)} -V_{n(k)}) 
(-\Delta +1)^{-1/2}\| =0  
\end{equation} 
\end{proposition} 

\begin{proof} Given $K$, let $P_K$ be the projection in momentum space onto $\abs{p}\leq 
K$ and $Q_K =1-P_K$. \eqref{4.12} implies that for any $\veps >0$, 
\begin{equation} \lb{4.14} 
\langle\varphi,\abs{V_n}\varphi\rangle \leq \veps \|\nabla\varphi\|^2 + C_\veps 
\|\varphi\|^2  
\end{equation} 
for a fixed $C_\veps$ and all $n$. This implies that 
\begin{equation} \lb{4.15} 
\| \abs{V_n}^{1/2} (-\Delta +1)^{-1/2} Q_K\|^2 \leq \veps +C_\veps (K^2 +1)^{-1/2}
\end{equation} 
so 
\begin{equation} \lb{4.16} 
\lim_{K\to\infty}\, \sup_n\, \|\abs{V_n}^{1/2} (-\Delta +1)^{-1/2} Q_K\| =0  
\end{equation} 

Thus, by a standard diagonalization argument, it suffices to show that for each $K$, there 
is a subsequence so that 
\begin{equation} \lb{4.17} 
\lim_{j,k\to\infty}\, \|(-\Delta+1)^{-1/2} P_K (V_{n(j)} -V_{n(k)}) P_K (-\Delta+1)^{-1/2}\| 
=0  
\end{equation} 

In momentum space, 
\begin{equation} \lb{4.18} 
Q_n = (-\Delta +1)^{-1/2} P_K V_n P_K (-\Delta +1)^{-1/2}  
\end{equation} 
has an integral kernel
\begin{equation} \lb{4.19} 
Q_n(p,q) =\chi_{\abs{p}\leq K} (p) (p^2 +1)^{-1/2} \widehat V_n (p-q) 
\chi_{\abs{q}\leq K} (q)(q^2 +1)^{1/2}  
\end{equation} 

By \eqref{4.12} and the fixed support hypothesis, we have 
\begin{equation} \lb{4.20} 
\sup_n\, (\|V_n\|_{L_1} + \|\vec x\, V_n\|_{L^1}) <\infty  
\end{equation} 
so that 
\begin{equation} \lb{4.21} 
\sup_n \, (\abs{\widehat V_n(k)} + \abs{\nabla\widehat V_n (k)}) <\infty  
\end{equation} 
which means $\{V_n(k) \mid \abs{k} \leq 2K\}$ is a uniformly equicontinuous family, 
so we can find a subsequence so 
\begin{equation} \lb{4.22} 
\lim_{j,k\to\infty}\, \sup_{\abs{k}\leq 2K}\, \abs{\widehat V_{n(j)}(k) - 
\widehat V_{n(\ell)}(k)} =0
\end{equation} 
It follows from \eqref{4.19} that 
\begin{equation} \lb{4.23} 
\int \abs{Q_{n(j)} (p,q) - Q_{n(\ell)} (p,q)}^2 \, dpdq\to 0  
\end{equation} 
so \eqref{4.17} holds since the Hilbert-Schmidt norm dominates the operator norm. 
\end{proof}

Given $V\in K_\nu$, we say $\wti H$ is a limit point at infinity in direction $e$ if 
there exists $x_n\to\infty$ with $x_n/\abs{x_n}\to e$ so that for the characteristic 
function of any ball and $z\in\bbC\backslash [a,\infty)$, we have 
\begin{equation} \lb{4.24} 
\lim_{n\to\infty}\, \|[(-\Delta +V(x-x_n)-z)^{-1}   -(\wti H-z)^{-1}]\chi_S\|=0
\end{equation} 
Let $\calL_e$ denote the set of limit points in direction $e$. Then our standard 
argument using Theorem~\ref{T4.2} and Proposition~\ref{P4.4} to get compactness 
implies 

\begin{theorem}\lb{T4.5} Let $V\in K_\nu$. Then 
\begin{equation} \lb{4.25} 
\sigma_\ess (-\Delta +V) = \ol{\bigcup_e\, \bigcup_{\wti H\in\calL_e} \sigma(\wti H)} 
\end{equation} 
\end{theorem}

\section{Potentials Asymptotic to Isospectral Tori} \lb{s5}

As a warmup, we will prove the following result which includes Theorem~\ref{T1.2} 
as a special case. We will consider functions $f\colon\bbR^\nu\to\bbR^\nu$ so 
\begin{equation} \lb{5.1} 
\lim_{\abs{x}\to\infty}\, \sup_{\abs{y}\leq L}\, \abs{f(x)-f(x+y)} =0  
\end{equation} 
for each $L$. For example, if $f$ is $C^1$ outside some ball and $\abs{\nabla f(x)} 
\to 0$ (e.g., $f(x) =\sqrt{x}\f{x}{\abs{x}}$), then \eqref{5.1} holds. 

\begin{theorem}\lb{T5.1} Let $V$ be a function on $\bbR^\nu$,
periodic in $\nu$ independent directions,
so $V$ is uniformly Kato {\rm{(}}e.g., $V\in L_\loc^p$ with $p$ a canonical 
value for $\bbR^\nu${\rm{)}}. Suppose either $V$ is bounded or $\abs{f(x)-f(y)} 
\leq (1-\veps)\abs{x-y}$ for some $\veps >0$. Let $f$ obey \eqref{5.1}. Let 
$W(x)=V(x+f(x))$. Then 
\begin{equation} \lb{5.2} 
\sigma_\ess (-\Delta +W) =\sigma (-\Delta +V)  
\end{equation} 
\end{theorem} 

\begin{remark} The condition that $V$ is bounded or $f$ is globally Lifschitz is 
needed to assure $W$ is locally $L^1$. We thank V.~Georgescu for pointing out to 
us the need for this condition, which was missing in our original preprint.  
\end{remark}

\begin{proof} Let $L$ be the integral lattice generated by some set of periods so 
$V(x+\ell)=V(x)$ if $\ell\in L$. Let $\pi\colon\bbR^\nu \to\bbR^\nu/L$ be the 
canonical projection. If $x_j\in\bbR^\nu$, since $\bbR^\nu/L$ is compact, we can 
find a subsequence $m(j)$ so $\pi((x_{m(j)}) + f(x_{m(j)}))\to x_\infty$. Then 
\[
-\Delta +W (\bddot -x_{m(j)}) \to -\Delta+V(x-x_\infty) 
\] 
so the limits are translates of $-\Delta +V$\!, which all have the same essential 
spectrum. \eqref{5.2} is immediate from Theorem~\ref{T3.11}. 
\end{proof} 

Our next result includes Theorem~\ref{T1.3}. 

\begin{theorem}\lb{T5.2} Let $W\colon\bbR^d\to\bbR$ be bounded and continuous, and 
obey 
\begin{equation} \lb{5.3} 
W(x+a)=W(x)  
\end{equation} 
if $a\in\bbZ^d$. Let $(\alpha_1, \dots, \alpha_d)$ be such that $\{(\alpha_1 n, 
\alpha_2 n,\dots, \alpha_d n)\mid n\in\bbZ\}$ is dense in $\bbR^d/\bbZ^d$ 
{\rm{(}}i.e., $1,\alpha_1,\dots,\alpha_d$ are rationally independent{\rm{)}}. 
Let $f\colon\bbZ\to\bbR^d$ obey 
\[
\lim_{n\to\infty}\, \sup_{\abs{m}\leq L}\, \abs{f(n) -f(n+m)} =0 
\] 
for each $L$. Let $V_0(n) = W(\alpha n)$ and let  
\begin{equation} \lb{5.4} 
V(n) =W(\alpha n + f(n))  
\end{equation} 
On $\ell^2(\bbZ)$, let $(h_0 u)(n)=u(n+1) + u(n-1)$. Then 
\begin{equation} \lb{5.5} 
\sigma_\ess (h_0 +V) =\sigma(h_0 +V_0)  
\end{equation} 
\end{theorem} 

\begin{proof} For each $x\in\bbR^d/\bbZ^d$, define 
\begin{equation} \lb{5.6} 
V_x(n) =W(\alpha n+x)  
\end{equation} 
Then a theorem of Avron-Simon \cite{AvS} (see \cite{CFKS}) shows that $\sigma (h_0+V_x)$ 
is independent of $x$ (and purely essential). Given any sequence $n_j$, find a 
sequence $n_{j(m)}$ so $f(n_{j(m)})\to x_\infty$ in $\bbR^d/\bbZ^d$. Then 
$V(n+n_{j(m)})\to V_{x_\infty}(n)$, so by Theorem~\ref{T1.6}, 
\[
\sigma_\ess (h_0 +V) = \ol{\bigcup_x \sigma (h_0 +V_x)=\sigma(h_0 +V_0)} 
\qedhere 
\]
\end{proof} 

Next, we turn to Theorem~\ref{T1.4} in the OPUC case. Any set of periodic Verblunsky 
coefficients $\{\alpha_n\}_{n=0}^\infty$ with 
\begin{equation} \lb{5.7} 
\alpha_{n+p} =\alpha_n
\end{equation} 
for some $p$ defines a natural function on $\bbC\backslash\{0\}$, $\Delta(z) =
z^{-p/2} \tr (T_p(z))$, where $T_p$ is a transfer matrix; see Section~11.1 of \cite{OPUC2}. 
(If $p$ is odd, $\Delta$ is double-valued; see Chapter~11 of \cite{OPUC2} for how to 
handle odd $p$.) $\Delta$ is real on $\partial\bbD$ and $\sigma_\ess (\calC(\alpha))$ 
is a union of $\ell$ disjoint intervals; $\ell\leq p$ (generically, $\ell=p$). As proven 
in Chapter~11 of \cite{OPUC2}, 
\begin{equation} \lb{5.8} 
\{\beta\in\bbD^p\mid\Delta (z;\{\beta_{n\text{ mod } p}\}_{n=0}^\infty) = 
\Delta(z;\alpha)\}\equiv T_\alpha  
\end{equation}  
is an $\ell$-dimensional torus called the isospectral torus. Moreover, the two-sided 
CMV matrix, defined by requiring \eqref{5.8} for all $n\in\bbZ$, has 
\begin{equation} \lb{5.9a} 
\sigma(\ti\calC(\beta))=\sigma_\ess (\calC(\alpha))  
\end{equation} 
for any $\beta\in T_\alpha$. 

Given two sequences $\{\kappa_n\}_{n=0}^\infty$ and $\{\lambda_n\}_{n=0}^\infty$ in $\bbD^p$, 
define 
\begin{equation} \lb{5.9} 
d(\kappa,\lambda) \equiv\sum_{n=0}^\infty e^{-n} \abs{\kappa_n -\lambda_n}  
\end{equation} 
Convergence in $d$-norm is the same as sequential convergence. We define 
\[
d(\kappa, T_\alpha) =\inf_{\beta\in T_\alpha}\, d(\kappa, \beta) 
\]
A sequence $\gamma_n$ is called asymptotic to $T_\alpha$ if 
\begin{equation} \lb{5.10} 
\lim_{m\to\infty}\, d(\{\gamma_{n+m}\}_{n=0}^\infty, T_\alpha) = 0  
\end{equation} 
Then the OPUC case of Theorem~\ref{T1.4} (settling Conjecture~12.2.3 of \cite{OPUC2}) 
says 

\begin{theorem} \lb{T5.3} Let \eqref{5.10} hold. Then 
\begin{equation} \lb{5.11} 
\sigma_\ess (\calC(\{\gamma_n\}_{n=0}^\infty)) = \sigma_\ess (\calC(\{\alpha_n\}_{n=0}^\infty))  
\end{equation} 
\end{theorem} 

\begin{proof} The right limit points are a subset of $\{\ti\calC 
(\{\beta_{n\text{ mod }p}\}_{n=-\infty}^\infty)\mid\{\beta\}_{n=0}^{p-1}\in T_\alpha\}$, so 
by Theorem~\ref{T3.1} and \eqref{5.9a}, \eqref{5.11} holds. 
\end{proof} 

By the same argument using isospectral tori for periodic Jacobi matrices 
\cite{FlMcL,Krich1,Krich2,vMoer} and for Schr\"odinger operators 
\cite{DubMatNov,Levit,McvM}, one has 

\begin{theorem}\lb{T5.4} If $T$ is the isospectral torus of a given periodic Jacobi  
matrix, $\ti J$, and $J$ has Jacobi parameters obeying 
\begin{equation} \lb{5.12} 
\lim_{n\to\infty}\, \min_{\ti a, \ti b\in T}\, \sum_{\ell=1}^\infty \, 
[\abs{a_{n+\ell} - \ti a_\ell} + \abs{b_{n+\ell} -\ti b_\ell}] e^{-\ell} =0  
\end{equation}  
then 
\begin{equation} \lb{5.13} 
\sigma_\ess (J) =\sigma (\ti J)  
\end{equation} 
\end{theorem} 

\begin{theorem}\lb{T5.5} Let $T$ be the isospectral torus of a periodic potential, $V_0$, 
on $\bbR$ and $V$\! on $[0,\infty)$ in $K_1$ and 
\begin{equation} \lb{5.14} 
\lim_{\abs{x}\to\infty}\, \inf_{W\in T}\, \int_0^\infty 
\abs{V (y+x)-W(y)} e^{-\abs{y}}\, dy =0  
\end{equation} 
then 
\begin{equation} \lb{5.15} 
\sigma_\ess \biggl( -\f{d^2}{dx^2} + V\biggr) = \sigma \biggl( -\f{d^2}{dx^2} +V_0\biggr)  
\end{equation} 
where $-\f{d^2}{dx^2} +V$ is defined on $L^2 (0,\infty)$ with $u(0)=0$ boundary conditions 
and $-\f{d^2}{dx^2}+V_0$ is defined on $L^2 (\bbR,dx)$. 
\end{theorem}

The following provides an alternate proof of Theorem~4.3.8 of \cite{OPUC1}: 

\begin{theorem} \lb{T5.6} Let $\{\alpha_j\}_{j=0}^\infty$ and $\{\beta_j\}_{j=0}^\infty$ 
be two sequences of Verblunsky coefficients. Suppose there exist $\lambda_j\in\partial\bbD$ 
so that 
\begin{alignat}{2}
\text{\rm{(i)}}& \qquad &\beta_j\lambda_j -\alpha_j &\to 0  \lb{5.17} \\ 
\text{\rm{(ii)}} &  \qquad & \lambda_{j+1} \bar\lambda_j &\to 1  \lb{5.18} 
\end{alignat}
Then 
\begin{equation} \lb{5.19} 
\sigma_\ess (\calC (\{\alpha_j\}_{j=0}^\infty)) = 
\sigma_\ess (\calC (\{\beta_j\}_{j=0}^\infty )) 
\end{equation} 
\end{theorem} 

\begin{proof} Let $\{\gamma_j\}_{j=-\infty}^\infty$ be a right limit of 
$\{\beta_j\}_{j=0}^\infty$, that is, $\beta_{\ell +n_k}\to\gamma_\ell$ for some 
$n_k$. By passing to a subsequence, we can suppose $\lambda_{n_j}\to\lambda_\infty$, 
in which case \eqref{5.18} implies $\lambda_{n_j+\ell} \to\lambda_\infty$ for 
each $\ell$ fixed. By \eqref{5.17}, $\{\lambda_\infty \gamma_j\}_{j=-\infty}^\infty$ 
is a right limit of $\{\alpha_j\}_{j=0}^\infty$. Since $\sigma (\ti\calC (\{\lambda 
\gamma_j\}_{j=-\infty}^\infty))$ is $\lambda$-independent, \eqref{5.19} follows from 
\eqref{3.12}. 
\end{proof}

\section{The HVZ Theorem} \lb{s6}

For simplicity of exposition, we begin with a case with an infinity-heavy particle; 
eventually we will consider a situation even more general than arbitrary $N$-body 
systems. Thus, $H$ acts on $L^2 (\bbR^{\mu(N-1)},dx)$ with 
\begin{equation} \lb{6.1} 
H=-\sum_{j=1}^{N-1} \, (2m_j)^{-1} \Delta_{x_j} + \sum_{j=1}^{N-1} V_{0j}(x_j) + 
\sum_{1\leq i < j \leq N-1} V_{ij}(x_j - x_i)  
\end{equation}  
where $x=(x_1, \dots, x_{N-1})$ with $x_j\in\bbR^\mu$. Here the $V$'s will 
be in $K_\mu$ with $K_\mu$ vanishing at infinity. $a$ will denote a partition 
$(C_1\dots C_\ell)$ of $\{0,\dots, N-1\}$ onto $\ell\geq 2$ clusters. We say $(ij)
\subset a$ if $i,j$ are in the same cluster, $C\in a$, and $(ij)\not\subset a$ if 
$i\in C_k$ and $j\in C_m$ with $k\neq m$, 
\begin{equation} \lb{6.2} 
H(a)=H-\sum_{\substack{ij\not\subset a \\ i <j}} V_{ij} (x_j -x_i)  
\end{equation} 
with $x_0\equiv 0$. The HVZ theorem says that 

\begin{theorem}\lb{T6.1} If each $V_{ij}$ is in $K_\mu$, $K_\mu$ vanishing at infinity, 
then 
\begin{equation} \lb{6.3} 
\sigma_\ess (H) = \ol{\bigcup_a \sigma (H(a))}  
\end{equation} 
\end{theorem} 

Since $H(a)$ commutes with translations of clusters, $H$ has the form $H(a)=T^a 
\otimes 1 + 1 \otimes H_a$ where $T^a$ is a Laplacian on $\bbR^{\mu(\ell-1)}$, 
and thus, if $\Sigma(a)=\inf\sigma (H_a)$, then $\sigma (H(a))= [\Sigma(a),
\infty)$. So \eqref{6.3} says 
\begin{equation} \lb{6.4} 
\sigma_\ess (H) = [ \Sigma,\infty) \qquad \Sigma\equiv\inf_a \,\Sigma(a) 
\end{equation}  

This result is, of course, well-known, going back to Hunziker \cite{Hun}, van Winter 
\cite{vW}, and Zhislin \cite{Zhi}, with geometric proofs by Enss \cite{Enss}, Simon \cite{S84}, 
Sigal \cite{Sig}, and Garding \cite{Gard}. Until Garding \cite{Gard}, all proofs involved some 
kind of combinatorial argument if only the existence of a Ruelle-Simon partition 
of unity. Like Garding \cite{Gard}, we will be totally geometric with a straightforward 
proof exploiting our general machine. $C^*$-algebra proofs can be found in Georgescu-Iftimovici 
\cite{GI1,GI3} and have a spirit close to our proof below. Rabinovich \cite{Rab05} 
has a proof of HVZ using his notion of invertibility at infinity that also has overlap 
with our philosophy. 

There is one subtlety to mention. Consider the case $\mu=1$, $N=3$, so
$\bbR^{\mu(N-1)} = \bbR^2 = \{(x_1,x_2)\mid x_1,x_2\in\bbR\}$.
There are then clearly six special directions: $\pm(1,0)$, $\pm(0,1)$, and 
$\pm(\f{1}{\sqrt2},\f{1}{\sqrt2})$. For any other direction $\hat e$, if $x_n/\abs{x_n} 
\to \hat e$, $V\to 0$, and the limit in that direction is $H_0 =H(\{0\},\{1\},\{2\})$. 

For $e=\pm (1,0)$, $\abs{(x_n)_1}\to\infty$ and $\abs{(x_n)_1 -(x_n)_2}\to\infty$, 
so the only limit at infinity would appear to be $H(\{0,2\},\{1\})$. But this is wrong! 
To say $x_n$ has limit $\pm(1,0)$ says $x_n/\abs{x_n}\to\pm (1,0)$, so $(x_n)_1\to\pm\infty$. 
But it does not say $(x_n)_2\to 0$, only $(x_n)_2/(x_n)_1\to 0$. For example, if 
$(x_n)_2\to\infty$, the limit is $H_0$. As we will see (it is obvious!), the limits are 
precisely $H_0$ and translates of $H(\{0,2\},\{1\})$. This still proves \eqref{6.3}, but 
with a tiny bit of extra thought needed. 

We want to note a general form for extending HVZ due to Agmon \cite{Agm}. We consider 
linear surjections $\pi_j\colon\bbR^\nu\to\bbR^{\mu_j}$ with $\mu_j\leq\nu$. Let 
$V_j\colon \bbR^{\mu_j}\to\bbR$ be in $K_{\mu_j}$ vanishing in $K_{\mu_j}$ sense 
at infinity. Then 
\begin{equation} \lb{6.4a} 
H=-\Delta+\sum_j V_j (\pi_j x)  
\end{equation} 
will be called an Agmon Hamiltonian. 

Given $e\in S^{\nu-1}$, define 
\begin{equation} \lb{6.5} 
H_e =-\Delta + \sum_{\{j\mid\pi_j e=0\}}\, V_j (\pi_j x)\equiv -\Delta +V_e  
\end{equation}  
Notice that since $H_e$ commutes with $x\to x+\lambda e$, $H_e$ has the form $H_e = 
-\Delta_e\otimes 1 + 1\otimes (-\Delta_{e^\perp} +V_e)$, so $\sigma (H_e)=[\Sigma_e, 
\infty)$ with $\Sigma_e =\inf\spec (H_e)$. 

In general, if $\cap_j \ker(\pi_j)\neq\{0\}$, $H$ has some translation invariant 
degrees of freedom and can, and should, be reduced, but the HVZ theorem holds for 
the unreduced case (and also for the reduced case, since the reduced $H$ which 
acts on $\bbR^\nu/\cap_j\ker(\pi_j)$ has the form \eqref{6.4a}). So we will not 
consider reduction in detail. 

By using $\pi_j$ to write $V_{ij} (x_i-x_j)$ in terms of mass scaled reduced coordinates, 
any $N$-body Hamiltonian has the form \eqref{6.4a}, and \eqref{6.4a} allows many-body forces.  
For the case of Theorem~\ref{T6.1}, if $e$ is given, define $a$ to be the partition 
with $(ij)\subset a$ if and only if $e_i=e_j$ (with $e_0\equiv 0$). Then $H_e =H(a)$ 
and \eqref{6.6} below is \eqref{6.3}. 

\begin{theorem}\lb{T6.2} For any Agmon Hamiltonian, 
\begin{equation} \lb{6.6} 
\sigma_\ess(H)= \ol{\bigcup_{e\in S^{\nu-1}} \sigma(H_e)}   
\end{equation} 
\end{theorem} 

\begin{proof} If $x_n/\abs{x_n}\to e$, we can pass to a subsequence where each $\pi_j 
x_n$ has a finite limit, or else has $\abs{\pi_j x_n}\to\infty$. It follows that the  
limit at infinity for $x_n$ is a translation (by $\lim \pi_j x_n$) of $H_e$ or of 
a limit at infinity of $H_e$. Thus, for any $\wti H$ in $\calL_e$,
the set of limits in direction $e$, 
\[
\sigma (\wti H)\subset \sigma(H_e) 
\]
and so, 
\[
\ol{\bigcup_{\wti H\in\calL_e}\, \sigma (\wti H})=\sigma (H_e) 
\]
and \eqref{6.6} is \eqref{4.25}. 
\end{proof} 

\begin{remark} It is not hard to see that as $e$ runs through $S^{\nu-1}$,
$\sigma(H_e)$ has only finitely many distinct values, so the closure in
\eqref{6.6} is superfluous.
\end{remark}

Because we control $\sigma_\ess (H)$ directly and do not rely on the a priori fact 
that one only has to properly locate $\inf\sigma_\ess(H)$ (as do all the proofs 
quoted above, except the original H,V,Z proofs and Simon \cite{S84}), we can obtain 
results on $N$-body interactions where the particles move in a fixed background 
periodic potential with gaps that can produce gaps in $\sigma_\ess (H)$.

\section{Additional Applications} \lb{s7}

We want to consider some additional applications of our machinery that shed light 
on earlier works: 
\begin{SL} 
\item[(a)] Sparse bumps, already considered by Klaus \cite{Klaus} using Birman-Schwinger 
techniques, by Cycon et al.\ \cite{CFKS} using geometric methods, and by Hundertmark-Kirsch 
\cite{HK} using methods that are essentially the same as the specialization of our 
argument to this example. Georgescu-Iftimovici \cite{GI3} also have a discussion of 
sparse potentials that overlaps our discussion. 

\item[(b)] Jacobi matrices with $a_n\to 0$ and CMV matrices with $\abs{\alpha_n}\to 1$ 
already studied by Maki \cite{Maki}, Chihara \cite{Chi70} (Jacobi), and by Golinskii 
\cite{Gol2000} (CMV). 

\item[(c)] Bounded Jacobi matrices and CMV matrices with finite essential spectrum 
already studied by Krein (in \cite{AK}) and Chihara \cite{Chi} (Jacobi case), and 
by Golinskii \cite{Gol2000} (CMV case). 
\end{SL} 

\begin{remark} Golinskii \cite{Gol2000} for (b) and (c) did not explicitly use CMV 
matrices  but rather studied measures on $\partial\bbD$, but his results are equivalent 
to statements about CMV matrices. 
\end{remark} 

Here is the sparse potentials result: 

\begin{theorem} [\cite{Klaus,CFKS}] \lb{T7.1} Let $W$\! be an $L^1$ potential of compact 
support on $\bbR$. Let $x_0 <x_1 <\cdots < x_n <\cdots$ so $x_{n+1}-x_n\to\infty$. Let 
\begin{equation} \lb{7.1} 
V(x)=\sum_{j=0}^\infty W(x-x_j)  
\end{equation} 
Then 
\begin{equation} \lb{7.2} \sigma_\ess \biggl( -\f{d^2}{dx^2}+V(x)\biggr) = 
\sigma\biggl( -\f{d^2}{dx^2}+W\biggr)  
\end{equation}  
\end{theorem} 

\begin{remarks} 1. That $W$\! has compact support is not needed. $W(x)\to 0$ sufficiently 
fast (e.g., bounded by $x^{-1-\veps}$) will do with no change in proof. 

2. Discrete eigenvalues of $-\f{d^2}{dx^2}+W$\! are limit points of eigenvalues for 
$-\f{d^2}{dx^2}+V$\!. 

3. There is a higher-dimensional version of this argument; see \cite{HK}. 
\end{remarks} 

\begin{proof} The limits at infinity are $-\f{d^2}{dx^2}$ and $-\f{d^2}{dx^2} +W(x-a)$. 
Now use Theorem~\ref{T3.11} or Theorem~\ref{T4.5}. 
\end{proof} 

\begin{remark} This example is important because it shows that one needs $\sigma 
(\wti H)$ and not just $\sigma_\ess (\wti H)$. 
\end{remark} 

As for $a_n\to 0$: 

\begin{theorem} [\cite{Chi70}]\lb{T7.2} Let $J$ be a bounded Jacobi matrix with $a_n\to 0$. 
Let $S$ be the limit points of $\{b_n\}_{n=1}^\infty$. Then 
\begin{equation} \lb{7.3} 
\sigma_\ess (J)=S  
\end{equation} 
\end{theorem} 

\begin{proof} The limit points at infinity are diagonal matrices with diagonal matrix 
elements in $S$, and by a compactness argument, every $s\in S$ is a diagonal matrix element 
of some limit. Theorem~\ref{T3.3} implies \eqref{7.3}. 
\end{proof} 

\begin{theorem} [\cite{Gol2000}]\lb{T7.3} Let $\calC(\{\alpha_n\}_{n=0}^\infty)$ be a CMV 
matrix of a sequence of Verblunsky coefficients with 
\begin{equation} \lb{7.3a} 
\lim_{n\to\infty}\, \abs{\alpha_n}=1  
\end{equation} 
Let $S$ be the set of limit points of $\{-\bar\alpha_{j+1}\alpha_j\}$. Then 
\begin{equation} \lb{7.4} 
\sigma_\ess (\calC(\{\alpha_j\}_{j=1}^\infty))=S  
\end{equation}  
\end{theorem} 

\begin{proof} By compactness of $\partial\bbD$, if $s\in S$, there is a sequence $n_j$ 
so $\alpha_{n_j+\ell}$ has a limit, $\beta_\ell$, for all $\ell$ and $s =-\bar\beta_1 
\beta_0$. The limiting CMV matrices have $\abs{\beta_\ell}=1$ by \eqref{7.3a}, so are 
diagonal with matrix elements $-\bar\beta_{\ell+1}\beta_\ell$. Thus, the spectra of 
limits lie in $S$, and by the first sentence, any such $s\in S$ is in the spectrum of 
a limit. Now use Theorem~\ref{T3.2}. 
\end{proof} 

Finally, we turn to the case of finite essential spectrum, first for Jacobi matrices. 

\begin{theorem} \lb{T7.4} Let $x_1, \dots, x_\ell\in\bbR$ be distinct. A bounded 
Jacobi matrix $J$ has 
\begin{equation} \lb{7.5} 
\sigma_\ess(J)=\{x_1, \dots, x_\ell\}  
\end{equation} 
if and only if 
\begin{SL} 
\item[{\rm{(i)}}] 
\begin{equation} \lb{7.6} 
\lim_{n\to\infty}\, a_n a_{n+1} \dots a_{n+\ell-1} =0  
\end{equation} 
\item[{\rm{(ii)}}] If $k\leq l$ and $n_j$ is such that 
\begin{equation} \lb{7.7}
 a_{n_j}\to 0 \qquad a_{n_j+k}\to 0 
\end{equation}
\begin{alignat}{2} 
& a_{n_j+m}\to\ti a_m\neq 0 \qquad\qquad && m=1,2,\dots, k-1 \lb{7.8} \\ 
& b_{n_j+m}\to \ti b_m \qquad\qquad && m=1,2,\dots, k \lb{7.9}  
\end{alignat} 
then the finite $k\times k$ matrix, 
\begin{equation} \lb{7.10} 
\ti J = \begin{pmatrix} 
\ti b_1 & \ti a_1 & {} \\ 
\ti a_1 & \ti b_2 & \ti a_2 \\
& \ddots & \ddots & \ddots  \\
& {} & \ddots & \ddots & \ti a_{k}  \\ 
& {} & {} & \ti a_{k-1} & \ti b_k
\end{pmatrix}
\end{equation} 
has spectrum a $k$-element subset of $\{x_1, \dots, x_\ell\}$. 

\item[{\rm{(iii)}}] Each $x_j$ occurs in at least one limit of the 
form \eqref{7.10} 
\end{SL} 
\end{theorem} 

\begin{proof} By Theorem~\ref{T3.3}, \eqref{7.5} holds if and only if the limiting 
$\ti J$'s have spectrum in $\{x_1, \dots, x_\ell\}$ and there is at least one $\ti J$ 
with each $x_j$ in the spectrum. $\ti J$ is a direct sum of finite and/or semi-infinite 
and/or infinite pieces. The semi-infinite pieces correspond to Jacobi matrices with 
nontrivial measures which have infinite spectrum. The two-sided infinite pieces also 
have infinite spectrum. Finite pieces of length $m$, which have $a$'s nonzero, have $m$ 
points in their spectrum, so no limit can have a direct summand of length $\ell+1$ or more.  
Thus, by compactness, \eqref{7.6} holds, that is, any set of $\ell$ $a$'s in the limit 
must have at least one zero. (ii) is then the assertion that the limits have spectrum 
in $\{x_1,\dots, x_\ell\}$, and (iii) is that each $x_j$ occurs. 
\end{proof} 

\begin{theorem}\lb{T7.5} 
\begin{SL} 
\item[{\rm{(a)}}] $J$ obeys  
\begin{equation} \lb{7.11a} 
\sigma_\ess(J) \subset \{x_1, \dots, x_\ell\}  
\end{equation} 
if and only if every right limit, $\ti J$, obeys 
\begin{equation} \lb{7.12} 
\prod_{j=1}^\ell (\ti J-x_j)\equiv P(\ti J)=0  
\end{equation} 

\item[{\rm{(b)}}] $J$ obeys \eqref{7.11a} if and only if $P(J)$ is compact. 
\end{SL} 
\end{theorem} 

\begin{proof} (a) \eqref{7.12} holds if and only if $\sigma (\ti J)\subset 
\{x_1, \dots, x_\ell\}$, so this follows from Theorem~\ref{T3.3}. 

\smallskip 
(b) $P(J)$ has finite width. Thus, it is compact if and only if all matrix elements 
go to zero, which is true (by compactness of translates of $J$) if and only if 
\eqref{7.12} holds for all limits. 
\end{proof} 

We have now come full circle---for Theorem~\ref{T7.5}(b) is precisely Krein's criterion 
(stated in \cite{AK}), whose proof is immediate by the spectral mapping theorem and 
the analysis of the spectrum of compact selfadjoint operators. However, our Theorem~\ref{T7.4} 
gives an equivalent, but subtly distinct, way to look at the limits. To see this, consider 
the case $\ell=2$, that is, two limiting eigenvalues $x_1$ and $x_2$.

This has been computed by Chihara \cite{ChiNATO}, who found necessary and sufficient 
conditions for $\sigma_\ess (J)=\{x_1,x_2\}$ are (there is a typo in \cite{ChiNATO}, 
where we give $(b_n-x_1) (b_n-x_2)$ in \eqref{7.13}; he gives, after changing to our 
notation, $(b_n-x_1) (b_{n+1}-x_2)$): 
\begin{align} 
& \lim_{n\to\infty} \, (a_n^2 + a_{n-1}^2 + (b_n -x_1)(b_n-x_2)) =0 \lb{7.13} \\ 
& \lim_{n\to\infty}\, (a_n (b_n+b_{n+1} -x_1 - x_2)) =0 \lb{7.14} \\ 
& \lim_{n\to\infty} \, (a_n a_{n+1}) =0 \lb{7.15}   
\end{align} 
To see this from the point of view of $(J-x_1)(J-x_2)$, note that 
\begin{align} 
\langle\delta_n, (J-x_1)(J-x_2)\delta_n\rangle &= a_n^2 + a_{n-1}^2 + 
(b_n-x_1) (b_n-x_2)  \lb{7.16}  \\ 
\langle\delta_{n+1}, (J-x_1)(J-x_2)\delta_n\rangle &= a_n (b_n-x_2) + 
a_n (b_{n+1}-x_1) \lb{7.17} \\ 
\langle\delta_{n+2}, (J-x_1)(J-x_n)\delta_n\rangle &= a_n a_{n+1} \lb{7.18} 
\end{align} 

If we think in terms of limit points, we get a different-looking set of equations. 
Consider limits, $\ti J$. Of course, \eqref{7.15} is common 
\begin{equation} \lb{7.19} 
\ti a_n \ti a_{n+1} =0  
\end{equation} 
But the conditions on summands of $\ti J$ become 
\begin{gather} 
\ti a_n = \ti a_{n-1} =0 \Rightarrow \ti b_n = x_1 \quad \text{or} \quad \ti b_n =x_2 \lb{7.20} \\ 
\ti a_n\neq 0 \Rightarrow \ti b_{n+1} + \ti b_n =x_1 + x_2 \quad\text{and}\quad 
\ti b_n \ti b_{n+1} - \ti a_n^2 = x_1 x_2 \lb{7.21}
\end{gather} 
For \eqref{7.20} is the result for $1\times 1$ blocks, and \eqref{7.21} says 
$2\times 2$ blocks have eigenvalues $x_1$ and $x_2$. It is an interesting exercise to 
see that \eqref{7.19}--\eqref{7.21} are equivalent to 
\begin{align}
& \ti a_n^2 + \ti a_{n+1}^2 + (\ti b_n -x_1) (\ti b_n-x_2) =0 \lb{7.22} \\ 
& \ti a_n (\ti b_n + \ti b_{n+1} -x_1 - x_2) =0 \lb{7.23} \\ 
& \ti a_n \ti a_{n+1} =0 \lb{7.24}  
\end{align} 

One can analyze CMV matrices similar to the above analysis. The analog of 
Theorem~\ref{T7.4} is: 

\begin{theorem}\lb{T7.6} Let $\lambda_1, \dots, \lambda_\ell\in\partial\bbD$ be distinct. 
A CMV matrix $\calC$ has 
\begin{equation} \lb{7.25} 
\sigma_\ess (\calC)=\{\lambda_1, \dots, \lambda_\ell\}  
\end{equation} 
if and only if 
\begin{SL} 
\item[{\rm{(i)}}] 
\begin{equation} \lb{7.26} 
\lim_{n\to\infty}\, \rho_n \rho_{n+1} \dots \rho_{n+\ell-1}=0  
\end{equation} 

\item[{\rm{(ii)}}] If $k\leq \ell$ and $n_j$ is such that 
\begin{alignat} {2}
& \rho_{n_j}\to 0 \qquad \rho_{n_j+k} \to 0  \lb{7.27} \\
& \alpha_{n_j+m}\to\ti\alpha_m \qquad\qquad && m=0,1,2,\dots, k-1, k \notag 
\end{alignat}
with $\abs{\ti\alpha_m}\neq 1$, $m=1, \dots, k-1$ {\rm{(}}by \eqref{7.27}, 
$\abs{\ti\alpha_0}=\abs{\ti\alpha_k}=1${\rm{)}}, then the matrix {\rm{(}}$\bdone = 
1\times 1$ unit matrix{\rm{)}}
\begin{equation} \lb{7.28} 
\ti\calC =[\Theta(\ti\alpha_1)\oplus\cdots\oplus\Theta(\ti\alpha_{k-1})] 
[-\ti\alpha_0 \bdone  \oplus\Theta(\ti\alpha_2)\oplus\cdots\oplus\Theta(\ti\alpha_{k-2}) 
\oplus\bar{\ti\alpha}_k\bdone]
\end{equation} 
if $k$ is even and 
\begin{equation} \lb{7.29} 
\ti\calC =[\Theta(\ti\alpha_1)\oplus\cdots\oplus\Theta(\ti\alpha_{k-2}) 
\oplus\bar{\ti\alpha}_k\bdone] [-\ti\alpha_0\bdone\oplus\Theta(\ti\alpha_2) 
\oplus\cdots\oplus\Theta (\ti\alpha_{k-1})] 
\end{equation} 
if $k$ is odd has eigenvalues $k$ elements among $\lambda_1, \dots,\lambda_\ell$. 

\item[{\rm{(iii)}}] Each of $\lambda_1, \dots, \lambda_\ell$ occurs as an eigenvalue 
of some $\ti\calC$. 
\end{SL} 
\end{theorem} 

\begin{proof} Same as Theorem~\ref{T7.4}. 
\end{proof} 

The analog of Theorem~\ref{T7.5} is 

\begin{theorem}\lb{T7.7} Let $\lambda_1, \dots, \lambda_\ell\in\partial\bbD$ be distinct. 
\begin{SL} 
\item[{\rm{(a)}}] $\calC$ obeys
\begin{equation} \lb{7.30a} 
\sigma_\ess (\calC)\subset \{\lambda_j,\dots,\lambda_\ell\}
\end{equation} 
if and only if every right limit $\ti\calC$ obeys 
\begin{equation} \lb{7.31} 
\prod_{j=1}^\ell (\ti\calC-\lambda_j) \equiv P(\ti\calC)=0 
\end{equation} 

\item[{\rm{(b)}}] $\calC$ obeys \eqref{7.30a} if and only if $P(\calC)$ is compact. 
\end{SL} 
\end{theorem} 

\begin{proof} Same as Theorem~\ref{T7.5}. 
\end{proof} 

We have come to Golinskii's OPUC analog of Krein's theorem \cite{Gol2000}. Again, 
it is illuminating to consider the case $\ell=2$. We will deal directly with 
limits of $\alpha_j$, call them $\ti\alpha_j$. The Theorem~\ref{T7.6} view of 
things is 
\begin{gather} 
\ti\rho_n \ti\rho_{n+1} =0 \lb{7.32} \\
\ti\rho_n =\ti\rho_{n+1}=0 \Rightarrow -\bar{\ti\alpha}_{n+1} \ti\alpha_n =
\lambda_1 \quad\text{or}\quad -\bar{\ti\alpha}_{n+1}\ti\alpha_n=\lambda_2 \lb{7.33} \\
\ti\rho_n\neq 0\Rightarrow -\bar{\ti\alpha}_n \ti\alpha_{n-1} - \bar{\ti\alpha}_{n+1} 
\ti\alpha_n=\lambda_1 + \lambda_2\quad\text{and}\quad \ti\alpha_{n-1} 
\bar{\ti\alpha}_{n+1} =\lambda_1\lambda_2 \lb{7.34} 
\end{gather}
\eqref{7.34} comes from the fact that the matrix $\calC$ of \eqref{7.28} is 
\begin{equation} \lb{7.35x} 
\begin{pmatrix} \bar{\ti\alpha}_n & \ti \rho_n \\ \rho_n & -\ti\alpha_n \end{pmatrix} 
\begin{pmatrix} -\ti\alpha_{n-1} & 0 \\ 0 & \bar{\ti\alpha}_{n+1}\end{pmatrix} 
\end{equation} 
where the determinant is $\ti\alpha_{n-1} \bar{\ti\alpha}_{n+1}$ and the trace is 
$-\bar{\ti\alpha}_n \ti\alpha_{n-1} - \ti\alpha_n \bar{\ti\alpha}_{n+1}$.  

  From the point of view of Theorem~\ref{T7.7}, using the CMV matrix is complicated 
since $(\calC-\lambda_1)(\calC-\lambda_2)$ is, in general, 9-diagonal! As noted by 
Golinskii \cite{Gol2000}, it is easier to use the GGT matrix (see Section~4.1 of 
\cite{OPUC1}), since it immediately implies 
\begin{equation} \lb{7.35} 
\ti\rho_n \ti\rho_{n+1} =\langle\delta_{n+2}, \calG^2 \delta_n\rangle =0 
\end{equation} 
and once that holds, $\calG$ becomes tridiagonal! Thus, one gets from $\langle 
\delta_{n+1}, (\calG -\lambda_1)(\calG-\lambda_2)\delta_n\rangle=0$ that 
\begin{equation} \lb{7.36} 
\ti\rho_n (-\bar{\ti\alpha}_n \ti\alpha_{n-1} - \bar{\ti\alpha}_{n+1} \ti\alpha_n 
-\lambda_1 -\lambda_2)=0 
\end{equation} 
and from $\langle\delta_n, (\calG-\lambda_1)(\calG-\lambda_2)\delta_n\rangle =0$, 
\begin{equation} \lb{7.37} 
(-\bar{\ti\alpha}_{n+1} \ti\alpha_n -\lambda_1)(-\bar{\ti\alpha}_n \ti\alpha_{n-1} 
-\lambda_2) - \ti\rho_n^2 \bar{\ti\alpha}_{n+1} \ti\alpha_{n-1} -\rho_{n+1}^2 
\bar{\ti\alpha}_{n-2} \ti\alpha_{n+1} =0 
\end{equation} 
Again, it is an interesting exercise that \eqref{7.32}--\eqref{7.34} are equivalent to 
\eqref{7.35}--\eqref{7.37}.

\section{Magnetic Fields} \lb{s8}

A magnetic Hamiltonian acts on $\bbR^\nu$ via 
\begin{equation} \lb{8.1} 
H(a,V)=-\sum_{j=1}^\nu \, (\partial_j -ia_j)^2 +V  
\end{equation}  
where $a$ is vector-valued. The magnetic field is the two-form defined by 
\begin{equation} \lb{8.2} 
B_{jk} =\partial_j a_k -\partial_k a_j  
\end{equation} 
If $\lambda$ is a scalar function, then 
\begin{equation} \lb{8.3} 
\ti a = a+\nabla\lambda  
\end{equation} 
produces the same $B$, and one has gauge covariance 
\begin{equation} \lb{8.4} 
H(\ti a,V)=e^{i\lambda}  H(a,V)e^{-i\lambda}
\end{equation} 
While the mathematically ``natural" conditions on $a$ are either $a\in L_\loc^4$, 
$\nabla \,\bddot\, a\in L_\loc^2$, or $a\in L_\loc^2$ (see \cite{CFKS,LeiS,S98}), for simplicity, 
we will suppose here that $B$ is bounded and uniformly H\"older continuous, that is, 
for some $\delta >0$, 
\begin{equation} \lb{8.5x} 
\sup_{x,j,k}\, \abs{B_{jk}(x)}<\infty \qquad 
\sup_{j,k,\abs{x-y}\leq 1}\, \abs{x-y}^{-\delta} \abs{B_{jk}(x)-B_{jk}(y)}<\infty  
\end{equation}  
It is certainly true that one can allow suitable local singularities. We will see 
later what \eqref{8.5x} implies about choices of $a$. With this kind of regularity on 
$B$, it is easy to prove that for a shift between different gauges of the type we 
consider below, the formal gauge covariance \eqref{8.4} is mathematically valid. 
Indeed, more singular gauge changes can be justified (see Leinfelder \cite{Lei}). 

If $a_j\to 0$ at infinity, it is easy to implement the ideas of Sections~\ref{s3} and 
\ref{snew4} with no change in the meaning of limit point at infinity; the limits all 
have no magnetic field. But as is well known, $a_j\to 0$ requires, very roughly 
speaking, that $B$ goes to zero at least as fast as $\abs{x}^{-1-\veps}$, so this 
does not even capture all situations where $B_{ij}\to 0$ at infinity. Miller \cite{Mill} 
(see also \cite{CFKS,S130}) noted that, in two and three dimensions, the way to control 
$B\to\ 0$ at infinity is to make suitable gauge changes in Weyl sequences---and that 
will also be the key to what we do here. 

We will settle for stating a very general limit theorem and not attempt to apply this 
theorem to recover the rather extensive literature on HVZ theorems and on essential 
spectra in periodic magnetic fields 
\cite{AHS,BrC,Cor,Hel,HM,HH,Hoe,Ift,IS,Iwa,Nak,Pas,UN,Vug,VZ93,VZ97,Zhi96a,Zhi96b,Zhi99,ZV}. 
We have no doubt that can be done and that the ideas below will be useful in 
future studies. We note that it should be possible to extend Theorem~\ref{T5.1} 
with ``slipped periodic" magnetic fields. 

\begin{definition} A set of gauges, $a_x$, depending on $x$ is said to be ``regular 
at infinity" if and only if, for every $R$, we have for some $\delta >0$, 
\begin{equation} \lb{8.5} 
\sup_{\abs{x-y}\leq R}\, \abs{a_x(y)} <\infty \qquad 
\sup_{\substack{ x,y,z \\ \abs{y-z}<1 \\ \abs{x-y} < R }}\, 
\abs{y-z}^{-\delta} \abs{a_x (y) - a_x(z)} < \infty
\end{equation} 
\end{definition} 

\begin{proposition} \lb{P8.1} If \eqref{8.5x} holds, there exists a set of gauges 
regular at infinity. 
\end{proposition} 

\begin{proof} The transverse gauge, $\vec a_{x_0}$, based at $x_0$ is defined by 
\begin{equation} \lb{8.6} 
a_{x_0;j} (x_0+y) = \sum_k \biggl[ \int_0^1 s B_{kj} (x_0+sy)\, ds \biggr] y_k
\end{equation} 
That this is a gauge is known (see below), and clearly, if $\abs{x_0-y} 
\leq R$, 
\[
\abs{\vec a_{x_0}(y)} \leq \tfrac12\, R\, \sup_x \, \|B(x)\| 
\]
and if $\abs{y-z}<1$ and $\abs{x_0-y}<R$, 
\[
\abs{\vec a_{x_0}(y) -\vec a_{x_0}(z)} \leq \tfrac12\, \bigl\{\sup_x\, \|B(x)\| + 
(R-1) \sup_{\abs{y-z}\leq 1}\, [\abs{y-z}^{-\delta} \|B(y) - B(z)\|]\bigr\} 
\qedhere 
\]
\end{proof} 

\begin{remarks} 1. We will call the choice \eqref{8.6} the local transverse gauge. 

2. Transverse gauge goes back at least to Uhlenbeck \cite{Uhl}, who calls them 
exponential gauge. They have been used extensively by Loss-Thaller \cite{LT} 
(see also Thaller \cite{Thal}) to study scattering.  

3. To see that \eqref{8.6} is a gauge is a messy calculation if done directly, 
but there is a lovely indirect argument of Uhlenbeck \cite{Uhl}. Without loss, 
take $x_0=0$. Call a gauge transverse if $\vec a(0)=0$ and $\vec x \,\bddot\, 
\vec a =0$. Transverse gauges exist, for if $\vec a_0$ is any gauge and 
\begin{equation} \lb{8.7a} 
\varphi(\vec x) = -\int_0^1 \vec x\, \bddot \, a_0 (s\vec x)\, ds
\end{equation}
then $\vec x \,\bddot\, \nabla\varphi = r \f{\partial}{\partial r} \varphi 
= -\vec x \,\bddot\, a_0(x)$, so $a=a_0 + \nabla\varphi$ is transverse. 
Next, note that if $\vec a$ is a transverse gauge, then 
\begin{align} 
\sum x_k B_{kj} &= (x\,\bddot\nabla) a_j - \vec\nabla_{\!j} (x\,\bddot\, a) + a_j \notag \\ 
&= \f{\partial}{\partial r}\, ra_j \lb{8.7b}
\end{align} 
Integrating \eqref{8.7b} shows \eqref{8.6} with $y=0$ is not only a gauge but the 
unique transverse gauge. 
\end{remarks}

If $a_x$ is a set of gauges regular at infinity, we say $\wti H$ is a limit at infinity 
of $H(a,V)$ in direction $\hat e$ if and only if with 
\begin{equation} \lb{8.7} 
(U_x\varphi)(y)=\varphi (y-x)  
\end{equation} 
we have that for some sequence $x_n$, $\abs{x_n}\to \infty$, $x_n/\abs{x_n}\to e$, and  
for each $R<\infty$ and $z\in\bbC\backslash [\alpha,\infty)$, 
\begin{equation} \lb{8.8} 
U_{x_n} ((H(a_{x_n},V)-z)^{-1}) U_{x_n}^{-1} \chi_R\to (\wti H-z)^{-1} \chi_R
\end{equation}  
with $\chi_R$ the characteristic function of a ball of radius $R$ about $0$. As usual, 
$\calL_e$ denotes the limits at infinity in direction $e$. 

\begin{theorem}\lb{T8.2} If $V\in K_\nu$ and $B$ obeys \eqref{8.5x}, then 
\begin{equation} \lb{8.9} 
\sigma_\ess   (H(a,V)) = \ol{\bigcup_{e\in S^{\nu -1}} \bigcup_{\wti H\in\calL_e} 
\sigma (\wti H)}
\end{equation} 
\end{theorem} 

In \eqref{8.9}, we get the same union if, instead of all regular gauges at infinity, we 
take only the local transverse gauges. 

\begin{proof} By using gauge-transformed Weyl sequences as in \cite{CFKS}, it is easy to 
see the right side of \eqref{8.9} is contained in $\sigma_\ess (H(a,V))$. To complete 
the proof, we need only show the right side, restricted to local transverse gauges,  
contains $\sigma_\ess (H(a,V))$.  

Localization extends effortlessly since $[j, H(a,V)]=\vec\nabla j\,\bddot\, (\vec\nabla 
- i\vec a) + (\vec\nabla -i\vec a)\,\bddot\, \vec\nabla j$ and $\|(\vec\nabla -i\vec a) 
\varphi\|^2$ is controlled by $H(a,V)$. Thus, we only need compactness of the 
gauge-transformed operators. Since \eqref{8.5} says the $a_x$'s translated to $0$ are 
uniformly equicontinuous, compactness of the $a$'s is immediate. $V$'s are handled as 
in Section~\ref{snew4}. 
\end{proof}

\bigskip

\end{document}